\documentclass[epic,11pt]{article}
\usepackage{amsfonts,amsmath,amscd,amssymb,amsthm}
\usepackage{fullpage}
\usepackage[alphabetic,initials]{amsrefs}
\usepackage[all]{xy}
\usepackage{hyperref} 
\usepackage{mathrsfs}
\input xy%
\usepackage{amssymb}

\def\convf{\hbox{\space \raise-2mm\hbox{$\textstyle      \bigotimes \atop \scriptstyle \omega$} \space}}
\def\0{{\bar 0}}
\def\1{{\bar 1}}
\def\C{{\mathbb C}}
\def\J{{\mathbb J}}
\def\I{{\mathbb I}}
\def\K{{\mathbb K}}
\def\Z{{\mathbb Z}}

\def\B{{\mathbb B}}

\def\F{{\mathbb F}}
\def\G{{\mathbb G}}

\def\pwg1{{\operatorname{PWG}}}
\def\pwg{{\operatorname{pwg}}}

\def\dist{{\operatorname{dist}}}

\newcommand{\ttk}{\mathtt{k}}
\newcommand{\ttr}{\mathtt{r}}
\newcommand{\tte}{\mathtt{e}}

\newcommand{\ttS}{\mathtt{S}}
\newcommand{\ttT}{\mathtt{T}}

\newcommand{\tl}{\triangleleft}

\newcommand{\itema}{\item[{{\rm(a)}}]}
\newcommand{\itemb}{\item[{{\rm(b)}}]}
\newcommand{\itemc}{\item[{{\rm(c)}}]}

\newcommand{\noi}{\noindent}
\newcommand{\ga}{\alpha}
\newcommand{\gb}{\beta}
\newcommand{\gc}{\gamma}

\newcommand{\Gd}{\Delta}
\newcommand{\Gt}{\Theta}
\newcommand{\gd}{\delta}

\newcommand{\gs}{\sigma}

\newcommand{\go}{\omega}

\newcommand{\gt}{\tau}

\newcommand{\gl}{\lambda}

\newcommand{\gr}{\rho}

\newcommand{\gep}{\epsilon}
\newcommand{\gth}{\theta}
\newcommand{\op}{\oplus}

\newcommand{\A}{\mathbb A}

\newcommand{\ot}{\otimes}

\newcommand{\fg}{\mathfrak{g}}\newcommand{\fgl}{\mathfrak{gl}}
\newcommand{\fsl}{\mathfrak{sl}}\newcommand{\osp}{\mathfrak{osp}}

\newcommand{\fh}{\mathfrak{h}}

\newcommand{\fb}{\mathfrak{b}}

\newcommand{\fn}{\mathfrak{n}}

\newcommand{\fp}{\mathfrak{p}}

\newcommand{\ff}{\footnote}
\newfont{\eufm}{eufm10 scaled\magstep1}

 \newcommand{\ti}{\times}

\newcommand{\cO}{\mathcal{O}}

\newcommand{\cD}{\mathcal{D}}

\newcommand{\cB}{\mathcal{B}}

\newcommand{\cS}{\mathcal{S}}

\newcommand{\cH}{\mathcal{H}}

\newcommand{\ey}{\end{eqnarray}}
\newcommand{\by}{\begin{eqnarray}}
\newcommand{\nn}{\nonumber}

\newcommand{\bco}{\begin{conjecture}}
\newcommand{\ba}{\begin{alg}}
\newcommand{\ea}{\end{alg}}
\newcommand{\eco}{\end{conjecture}}
\newcommand{\bpf}{\begin{proof}}
\newcommand{\epf}{\end{proof}}
\newcommand{\bt}{\begin{theorem}}

\newcommand{\et}{\end{theorem}}

\newcommand{\br}{\begin{rem}}
\newcommand{\er}{\end{rem}}
\newcommand{\brs}{\begin{rems}}
\newcommand{\ers}{\end{rems}}
\newcommand{\bi}{\begin{itemize}}
\newcommand{\ei}{\end{itemize}}
\newcommand{\bl}{\begin{lemma}}
\newcommand{\bsul}{\begin{sublemma}}
\newcommand{\esul}{\end{sublemma}}
\newcommand{\bp}{\begin{proposition}}
\newcommand{\be}{\begin{equation}}
\newcommand{\bc}{\begin{corollary}}
\newcommand{\bexs}{\begin{examples}}
\newcommand{\eexs}{\end{examples}}
\newcommand{\bexa}{\begin{example}}
\newcommand{\eexa}{\end{example}}
\newcommand{\bex}{\begin{exercise}}
\newcommand{\eex}{\end{exercise}}
\newcommand{\btab}{\begin{tab}}
\newcommand{\etab}{\end{tab}}
\newcommand{\bg}{\begin{fig}}
\newcommand{\eg}{\end{fig}}
\newcommand{\el}{\end{lemma}}
\newcommand{\ep}{\end{proposition}}
\newcommand{\ee}{\end{equation}}
\newcommand{\ec}{\end{corollary}}
\newcommand{\Bc}{\begin{center}}
\newcommand{\Ec}{\end{center}}
\newcommand{\bh}{\begin{hyp}}
\newcommand{\eh}{\end{hyp}}
\newcommand{\bhs}{\begin{hyps}}
\newcommand{\ehs}{\end{hyps}}
\newcommand{\bd}{\begin{dfn}}
\newcommand{\ed}{\end{dfn}}

\begin{document}
\title{Table of Contents}

\newtheorem{thm}{Theorem}[section]
\newtheorem{hyp}[thm]{Hypothesis}
 \newtheorem{hyps}[thm]{Hypotheses}
  \newtheorem{rems}[thm]{Remarks}

\newtheorem{conjecture}[thm]{Conjecture}
\newtheorem{theorem}[thm]{Theorem}
\newtheorem{theorem a}[thm]{Theorem A}
\newtheorem{example}[thm]{Example}
\newtheorem{examples}[thm]{Examples}
\newtheorem{corollary}[thm]{Corollary}
\newtheorem{rem}[thm]{Remark}
\newtheorem{lemma}[thm]{Lemma}
\newtheorem{sublemma}[thm]{Sublemma}
\newtheorem{cor}[thm]{Corollary}
\newtheorem{proposition}[thm]{Proposition}
\newtheorem{exs}[thm]{Examples}
\newtheorem{ex}[thm]{Example}
\newtheorem{exercise}[thm]{Exercise}
\numberwithin{equation}{section}%
\setcounter{part}{0}
\newcommand{\drar}{\rightarrow}
\newcommand{\lra}{\longrightarrow}
\newcommand{\rra}{\longleftarrow}
\newcommand{\dra}{\Rightarrow}
\newcommand{\dla}{\Leftarrow}
\newcommand{\lfa}{\leftrightarrow}

\newtheorem{Thm}{Main Theorem}


\newtheorem*{thm*}{Theorem}
\newtheorem{lem}[thm]{Lemma}
\newtheorem{fig}[thm]{Figure}
\newtheorem*{lem*}{Lemma}
\newtheorem*{prop*}{Proposition}
\newtheorem*{cor*}{Corollary}
\newtheorem{dfn}[thm]{Definition}
\newtheorem*{defn*}{Definition}
\newtheorem{notadefn}[thm]{Notation and Definition}
\newtheorem*{notadefn*}{Notation and Definition}
\newtheorem{nota}[thm]{Notation}
\newtheorem*{nota*}{Notation}
\newtheorem{note}[thm]{Remark}
\newtheorem*{note*}{Remark}
\newtheorem*{notes*}{Remarks}
\newtheorem{hypo}[thm]{Hypothesis}
\newtheorem*{ex*}{Example}
\newtheorem{prob}[thm]{Problems}
\newtheorem{conj}[thm]{Conjecture}

\title{Explicit expressions for \v Sapovalov elements in Type A.}
\author{Ian M. Musson\ff{Research partly supported by  NSA Grant H98230-12-1-0249, and Simons Foundation grant 318264.} \\Department of Mathematical Sciences\\
University of Wisconsin-Milwaukee\\ email: {\tt
musson@uwm.edu}}
\maketitle
\begin{abstract}
We give explicit expressions for \v Sapovalov elements in Type A Lie algebras and superalgebras. Explicit expressions were already given in \cite{M2} Section 9, using non-commutative determinants, and in fact our first main results, Theorems \ref{bb} and \ref{bsb1}  can be viewed as complete expansions of these determinants. But we give new  proofs, which seem easier because they avoid induction and cofactor expansion. 
We also  describe \v Sapovalov elements for $\fgl(m,n)$ with respect to an arbitrary Borel subalgebra in Theorem \ref{bsb9}, and interpret 
\v Sapovalov elements in Type A as determinants of Hessenberg matrices in Theorems  \ref{hea} and \ref{7t}.The exact form of the explicit expressions  depends on an ordering on the set of positive roots, and Hessenberg matrices  are useful in changing the ordering.  Having explicit expressions for \v Sapovalov elements allows us to give easy proofs of several results in representation theory, see Section \ref{1surv} and Subsection \ref{mrs}.
\end{abstract}
\ref{se1} {Definition and significance of \v Sapovalov  elements.}\\
\ref{1s.8} {The Type A Case.}\\
\ref{sab} {\v Sapovalov elements for arbitrary Borel subalgebras.}\\ 
\ref{shm} {\v Sapovalov  elements as determinants of Hessenberg matrices.}\\
\ref{pse} {Powers of \v Sapovalov elements.}\\ 
\ref{1surv} {{Survival of \v{S}}apovalov elements in factor modules.}
\section{Definition and significance of \v Sapovalov  elements.}\label{se1}

Let $\fg$ be a simple Lie algebra or a contragredient Lie superalgebra with set of simple and positive roots $\Pi$ and $\Delta^+$ respectively.  In the Type A case we consider the cases $\fgl(m)$ and $\fgl(m,n)$. Fix a positive root $\eta$  and a  positive   integer $  \underline{m}$. If $\eta$ is isotropic, assume $\underline{m}=1$, and if $\eta$ is odd non-isotropic, assume that $\underline{m}$ is odd.
Let ${\overline{\bf P}}(\underline{m}\eta)$ be the set of partitions of $\underline{m}\eta$, see
\cite{M} Remark 8.4.3, and also (8.4.4) for the notation  $e_{-\pi}$ below. Then let
$\pi^0 \in {\overline{\bf P}}(\underline{m}\eta)$ be the unique partition of $\underline{m}\eta$ such
that $\pi^0(\ga) = 0$ if $\ga \in \Delta^+ \backslash \Pi.$
We say that $\gth = \theta_{\eta,\underline{m}}\in U({\mathfrak b}^{-} )^{- \underline{m}\eta}$ is a
{\it \v Sapovalov element for the pair} $(\eta,\underline{m})$ if it has the form
\be \label{rat}
\theta = \sum_{\pi \in {\overline{\bf P}}(\underline{m}\eta)} e_{-\pi}
H_{\pi},\ee where 
\be \label{rbt} H_{\pi} \in U({\mathfrak h}), \; H_{\pi^0} = 1,\ee  and
\be \label{boo} e_{\ga} \theta \in U({\mathfrak g})(h_{\eta} + \rho(h_{\eta})-\underline{m}(\eta,\eta)/2)+U({\mathfrak g}){\mathfrak n}^+ , \; \rm{ for \; all }\;\ga \in \Delta^+. \ee
For a semisimple Lie algebra, the existence of such elements was shown by \v Sapovalov, \cite{Sh} Lemma 1.
Consider the hyperplane in ${\mathfrak h}^*$ given by 
 \be \label{vat}{\mathcal H}_{\eta, \underline{m}} = \{ \lambda \in  {\mathfrak h}^*|(\lambda + \rho, \eta) = \underline{m}(\eta, \eta)/2  \}.\ee
From \eqref{boo}, the \v Sapovalov element
$\theta_{\eta,\underline{m}}$
has the important property that if $\gl
\in {\mathcal H}_{\eta, \underline{m}}$ then
$\theta_{\eta, \underline{m}}v_{\gl}$  is a highest weight vector of weight $\lambda -\underline{m}\eta$ in
$M(\gl)$.
The normalization condition $H_{\pi^0}=1$
guarantees that $\theta_{\eta, \underline{m}}v_\gl$ is never zero.
Hence  $xv_{\gl-m\eta}\lra  x\theta_{\eta,\underline{m}}v_{\gl}$ defines  a non-zero map between Verma modules $M(\gl-\underline{m}\eta) \lra M(\gl)$.
In the semisimple case, all maps between Vermas are composites of these maps,
\cite{H2} Theorem 5.1.\\ \\
\v Sapovalov elements  have also appeared in a number of situations in representation theory, usually in Type A.
Though not given this name, they appear in the work  of Carter and Lusztig \cite{CL}.
Indeed determinants similar to those in our  Theorem \ref{hea} were introduced in \cite{CL} Equation (5), and   our Theorem   \ref{11.5} may be viewed as a version of \cite{CL} Theorem 2.7. Carter and Lusztig use their result to study tensor powers of the defining representation of $GL(V)$, and homomorphisms between Weyl modules in positive characteristic, see also \cite{CP} and \cite{F}. Later Carter \cite{Car} used \v Sapovalov
elements to construct raising and lowering operators for $\fsl(n,\C)$, see also \cite{Br3}, \cite{Carlin}.  In \cite{Car}, these operators are  used to construct orthogonal bases for non-integral Verma modules, and all finite dimensional modules for $\fsl(n,\C)$.
\\ \\
From the results in Section \ref{pse} it is only really necessary to consider the case $\underline{m}=1$, so to simplify notation we set $\mathcal{H}_{\eta} = \mathcal{H}_{\eta,1},$ and
denote a  \v Sapovalov element for the pair $(\eta, 1)$ by $\gth_{\eta}$. In addition in the Type A case, it is only necessary to consider the highest root (resp. highest odd root) 
 of $\fgl(m)$ (resp. $\fgl(m,n)$), since any positive root of the same type is the highest root of some general linear subalgebra in the Lie algebra (resp. superalgebra) case.  These two observations greatly simplify our computations. 
If $\fg=\fgl(m,n)$) and we use the distinguished Borel subalgebra, then 
a \v Sapovalov  element for an even root is just 
a \v Sapovalov  element for  one of the summands of $\fg_0$.  We  do not consider 
what happens to  \v Sapovalov  elements for an even root when the Borel subalgebra is changed, but see \cite{M2} Section 4.
\\ \\
  Let $\fh$ be a Cartan subalgebra of $\fg$ and $W$ the Weyl  group. Unless otherwise stated we use the   distinguished Borel subaalgebra of $\fg$ and denote the corresponding set of simple and positive roots by  $\Pi$ and $\Gd^+$ respectively.  The even and odd roots are  
$\Gd_0$ and $\Gd_1$, and we set $\Gd_i^+ =\Gd_i \cap\Gd^+$, ${\gr_i}=\frac{1}{2}\sum_{\ga \in \Gd^+_i} \ga$ for $i=0, 1$ and ${\gr}={\gr_0}-{\gr_1}$.  The dot action of $W$  on $\fh^*$ is defined by $w\cdot \gl = w(\gl +\gr)-\gr$ for $w\in W$  and  $\gl\in\fh^*$. 
\\ \\
  For a basic classical simple Lie superalgebra $\fg
$ not of type A, the behaviour of \v Sapovalov  elements is rather subtle, Let  $\Pi$ be the distinguished or anti-distinguished set of roots 
and let $W_{\rm nonisotropic}$ (resp. $W_{\rm even}$) be
the subgroup of $W$ generated by the reflections $s_\ga,$ where $\ga$ is simple and {\rm nonisotropic} (resp {\rm even}).  Then in \cite{M2} \v Sapovalov  elements 
$\theta_{\gamma,m}$ are shown to exist under suitable conditions on $m$ if $\gc$ is in the $W_{\rm nonisotropic}$ or  $W_{\rm even}$ orbit of a simple root.
However even for $\fg=\osp(3,2)$  the latter condition does not always hold.  In  
the orthosymplectic case, 
we can gain some insight into the structure of Verma modules by considering the effect of changing the Borel subalgebra on \v Sapovalov  elements, \cite{M2} Section 4.  But complications can arise.  Consider  $\osp(3,2),$ which seems to be a good test case for these kind of issues.  Changing the Borel subalgebra leads to a map  Verma modules $M(s_\gc\cdot \gl) \lra M(\gl) $ which is non-zero and not injective with $\gc$ an even root, \cite{M} Example 9.3.4. Verma modules for  $\osp(3,2)$ can have 8 non-isomorphic composition factors \cite{Mas}, and 20 or 23 submodules \cite{M9}. 
\\ \\
This paper replaces Sections 9 and 10 of the unpublished manuscript \cite{M2}.  Earlier sections will be replaced by \cite{M3} and \cite{M4}.
I thank Jon Brundan for  suggesting that \v Sapovalov  elements could be used to prove Theorem \ref{1shapel}. I am also grateful to the referee for many helpful comments. 
\section{The Type A Case.}\label{1s.8}
The main results of this Section are Theorems \ref{bb} and \ref{bsb1}.
Let $\ga_i=  \gep_i -\gep_{i+1}$  be the simple roots of $\fgl(m)$ for $1 \le i \le m-1$.
For the superalgebra case we will also need the simple roots
$\gc_i=  \gd_i -\gd_{i+1}$ of $0\op\fgl(n)$ for $1 \le i \le n-1$,
and $\gb=\gep_m-\gd_1$. Set 
\by \label{heb} e_{\ga_i}= e_{i,i+1},\;  e_\gb &=& e_{m,m+1},\; e_{\gc_j} = e_{m+j,m+j+1},\\
h_{\ga_i}= e_{i,i} - e_{i+1,i+1},\; h_\gb &=& e_{m,m}+ e_{m+1,m+1}
,\; 
h_{\gc_j} =  e_{m+j+1,m+j+1} -e_{m+j,m+j}.
\nn\ey For any simple root ${\gs}$, let $e_{-\gs}$ be the transpose of $e_{\gs}$.  As remarked above,  
we can assume that $\eta =\gep_1 -\gep_{m}$  in the  $\fgl(m)$ case, and $\eta =\gep_1 -\gd_{n}$ in the case of $\fgl(m,n)$. The general case follows from this by relabelling the indices.
Let $(\;,\;)=(\;,\;)_{m,n}$ be the bilinear form on $\fh^*$ defined by
\be \label{hcb} (\gep_i,\gep_j)_{m,n} = - (\gd_i,\gd_j)_{m,n} = \gd_{i,j}\ee
for all relevant indices $i, j$.

\br \label{5r}{\rm 
\bi \itema
For a non-isotropic root $\ga$, set $\ga^\vee = 2\ga/(\ga,\ga)$.  From \eqref{hcb}, it follows that 
 $(\ga,\ga) = 2$, and  $\ga^\vee = \ga$,  (resp. $(\gc,\gc) = -2$ and $\gc^\vee = -\gc$)  for a simple root $\ga$ of $\fgl(m)\op 0$ (resp. $0\op\fgl(n)$).  In the latter case, it is not possible to
choose elements $e = e_{\gc}, f = e_{-\gc}$ and
$h = h_{\gc} =[e_{\gc}, e_{-\gc}]$, such that the elements of $e, f, h$ satisfy the usual defining relations for $sl(2)$, see for example \cite{M} (A.4.2).
 With the above choices for $\gc$ we have $[e_{\gc}, e_{-\gc}] =-h_{\gc}$. \ei}
\er
\noi For $\ga \in \fh^*$, let $h_\ga \in \fh$ be the unique element such that $(\ga,\gb) = \gb(h_\ga)$ for all $\gb \in \fh^*$.
\\ \\
If $N$ is a positive integer, we set $[N] =\{1,2,\ldots, N\}$. For $k\in [m-1]$, let 
$\gs_k= \gep_1-\gep_{k+1},$ 
and
\be \label{rt4} h_k= h_{\gs_k} +(\gr,{\gs_k}) -1.\ee

\bl \label{flo}
\bi \itema If $\ga$ is a simple root of $\fg =\fgl(m,n)$ then 
$2(\gr,\ga) = (\ga,\ga).$
\itemb For $k\in [m-1]$, $h_k -h_{k-1}=h_{\ga_k}+ 1$ if $k>1$, and $h_1=h_{\ga_1}$.
\ei
\el
\bpf (a) is well known, see for example \cite{M} Corollary 8.5.4.  Then (b) follows from \eqref{rt4}.\epf
\noi Now 
consider a strictly descending sequence of integers
\be \label{dod} I = \{i_0 , i_1, \ldots i_s, i_{s+1} \}.\ee
If $I$ is a  singleton set, then set $f_I =1$. Otherwise define  $f_{I} \in U(\fn^-)$ by
 \by \label{ddd}f_I = e_{i_{0},i_1}e_{i_1,i_2} e_{i_2,i_3} \ldots e_{i_s,i_{s+1}}. \ey
\subsection{The case of  $\fgl(m)$.}
Let
\be \label{120} \I = \{I\subseteq [m]|1, m\in I\}, \ee  
and for $ I\in \I$, define
\be \label{130} r(I) = \{s-1|s \in \bar I\} \ee
where $\bar I$ be the complement of $I$ in $\I$. Recall $h_i$ from  \eqref{rt4} and  define
\be \label{1q} H_J =\prod_{i \in r(J)}h_{_i}.\ee

\bt\label{bb} Let $\eta =\gep_1 -\gep_{m}$, and suppose $v$ is a highest weight vector in the $\fg$-Verma module $M(\gl)$ of weight $\gl$, and set
\be \label{cat22} \Gt_\eta= \sum_{J\subseteq \I} f_J H_J.\ee
Then $e_{\ga_k}\Gt_\eta v=0$ for $k\in [m-2]$.
If $(\gl+\gr,\eta)=1,$  Then $\Gt_\eta v$ is a $\fg$-highest weight vector, and so
$\Gt_\eta$ is a \v Sapovalov element for the pair $(\eta,1)$.
\et
\bpf
Clearly  sets in  $\I$ correspond to partitions of $\eta$, and \eqref{rat} holds. Furthermore the partition $\pi^0$ corresponds to $[m] \in \I$ and $e_{\pi^0}= f_{[m]}$. Since  $H_{[m]} =1$, \eqref{rbt} holds. 
We show that $\sum_{J\subseteq \I} f_J H_J$ satisfies Equation \eqref{boo} in the definition of a \v Sapovalov element.
Given a simple root $\ga_k=\gep_{k} - \gep_{k+1}$, to avoid double subscripts we sometimes set $\ga= \ga_k$.
Thus  $e_\ga= e_{k,k+1}$.
The goal is to show that if $v $  
iis a highest weight vector with weight  $\gl\in\cH_\eta$, then  $e_\ga \sum_{J\subseteq \I} f_J H_J v=0$.
Consider the set
\[\cS_k=\{I\subseteq \I|k,k+1 \in I   \}.\]
\noi
Suppose first $I \in \cS_k$. In this case $e_{-\ga}$ is a factor of $f_I.$
There is a unique $g$ such that $k+1 = i_{g}$ and $k = i_{g+1}$.
With $I$ as in \eqref{dod}, $i_0=m$ and $i_{s+1}=1$, we set  
\be \label{ay}I_1 = \{m, i_1, \ldots i_{g-1},k+1 \}\mbox{ and } I_2 = \{k, i_{g+2}, \ldots i_s, 1 \}. \ee
Then  $I$ is a  disjoint union $I = I_1 \cup I_2$, and we have
\be \label{2vam}
 e_\ga f_{I}v = f_{I_1} h_\ga f_{I_2}v =
\left\{ \begin{array}
  {cc}f_{I_1}h_\ga v&\mbox{if} \;\; I_2 = \{1\} \mbox{ is a singleton }
\\ f_{I_1}f_{I_2} (h_\ga +1)v &\mbox{ otherwise}.
\end{array} \right. \ee
\noi
\noi Define 
\be\label{qrt} I^+ = I \backslash \{k\}, \mbox{ and } I^- = I \backslash \{k+1\}.\ee
If $k=m-1$, then $I^- \notin \I$ and if $k=1,$ then $I^+ \notin \I$. We set
$f_{I^-}=0$ or $f_{I^+}=0$ in these cases. It is easy to check the following result.
\bl \label{hog} We have \bi \itema $e_\ga f_{I^+}v = f_{I_1} f_{I_2}v$
and $e_\ga f_{I^-}v = -f_{I_1} f_{I_2}v$.
\itemb \[r(I^+)=r(I)\cup \{k-1\} \mbox{ and } r(I^-)=r(I)\cup \{k\}\]
\itemc If $J\subseteq \I$ and $e_\ga f_{J}v \in  f_{I_1} f_{I_2}U(\fh)v$, then
 $$J=I, I^+\mbox{ or } I^-.$$ If $I^+$ or  $I^-$ is not a subset of $\I$, they should be left out in the above equation.
\ei \el
\noi From (b) in the Lemma and the definition of $H_J$, we have
\be \label{dog22} H_{I^+} = h_{k-1}H_I \mbox{ and } H_{I^-} = h_{k}H_I.\ee
Now we apply $e_\ga$ to the sum of terms in \eqref{cat22} corresponding to $I, I^+$ and $I^-$. 
Suppose that $k \in [m-2]$.  
For any highest weight vector $v$, we have using \eqref{dog22}, \eqref{2vam} and Lemma \ref{hog},
\be \label{aa}
e_\ga( f_{I}H_I + f_{I^+}H_{I^+} + f_{I^-}H_{I^-})v=\left\{ \begin{array}
  {cc}  f_{I_1}f_{I_2} (h_{\ga_1} -h_{1} ) H_I v &\mbox{if} \;\; k =1
\\  f_{I_1}f_{I_2} (h_{\ga_k} +1+h_{k-1} -h_{k} ) H_I v &\mbox{if} \;\; k>1.
\end{array} \right. \ee
 and this is  zero by Lemma
\ref{flo},  provided that $k \neq m-1$.
Now if $J$ does not contain $k$ or $k+1$, then
$e_\ga  f_J H_J v=0$. Otherwise
there exists a unique set $I\in \cS_k$ 
 such that
$J = I, I^+$ or $I^- $. Thus we have using \eqref{aa},
\[ e_\ga \sum_{J\subseteq \I} f_J H_J v= e_\ga \sum_{I \in \cS_k} (f_{I}H_I + f_{I^+}H_I^+ +f_{I^-}H_I^-)v =0.\]
This proves the first statement in the Theorem. Now suppose $k = m-1$, so $f_{I^-}=0$,  and  $\gl\in\mathcal{H}_{\eta} $. 
For $I \in \cS_{m-1}$, we have $e_\ga f_I v = h_\ga f_J v = f_J (h_\ga+1)v$ where $J   = I \backslash \{m\}$. Set $I^+ = I \backslash \{m-1\}$. Then $e_\ga f_{I^+} v = f_J v$.
Now by Equation
\eqref{dog22}, $H_{I^+} = h_{m-2}H_I $. Thus using \eqref{rt4} and the facts that $
{\gs_{m-2}}$ is a sum of ${m-2}$ simple roots, and 
${\ga_{m-1}} + {\gs_{m-2}} = \eta$, we have
\by \label{as} e_\ga (f_{I}H_I + f_{I^+}H_{I^+})v &=&
 f_J H_I (h_{\ga_{m-1}}+1 + h_{m-2})v \nn\\
&=& f_J H_I (h_{\ga_{m-1}} + h_{\gs_{m-2}} + m-2)v\nn\\
&=& f_J H_I (h_\eta+{m-2})v.
\ey Now since $\gl\in \cH_\eta$ and $(\gr, \eta) = m-1$,  it follows that $(\gl, \eta) = 2-m$.  Thus the expression in \eqref{as} is zero.
\epf
\br {\rm There is a simple method that can be used to determine the linear factors $h_k$ of $H_I$. Consider the case that $I=[m]$ is as large as possible. Then $H_I=1.$ If $\ga= \ga_k$, for $k\in [m-2]$ and $k>1$, consider the condition,  compare \eqref{aa},
\by \label{as4} 0 &=&e_\ga( f_{I}H_I + f_{I^+}H_{I^+} + f_{I^-}H_{I^-})v \nn\\ &=&
 f_{I_1} f_{I_2}(h_\ga+1 + h_{k-1} -h_k)v.\nn\\
\ey   This gives a recurrence relation on the coefficients $h_k$. Furthermore when $k=1$, there are only two terms in the above equation, and this gives $h_1$. This remark can be applied to other situations, in particular we used it to develop the proof of Theorem  \ref{bsb9}, see Lemma \ref{by2}. The recurence in \eqref{as4} is obtained by working from left to right along the Dynkin diagram.  We could have instead worked from  right to left, and obtained an different expression $ \sum_{J\subseteq \I} f_J H'_J$ for the \v Sapovalov element $\gth_\eta$.
(Recall that by \cite{M2} Theorem 2.1, 
\v Sapovalov elements are only unique modulo a left ideal.) However for 
$\gl\in \cH_\eta,$
$ \sum_{J\subseteq \I} f_J H_J(\gl) = \sum_{J\subseteq \I} f_J H'_J(\gl).$ Thus both expressions give rise to the same highest weight vector in 
$M(\gl) $. In the super case, Theorem \ref{bsb1} we work from both ends of the diagram towards the unique grey node.  In all the proofs the assumption 
$\gl\in \cH_\eta$ is used only in the last step.}\er

\subsection{The case of  $\fgl(m,n)$.} \label{cmn}
For $i\in [n-1]$, set $\gt_{i} = \gd_i - \gd_n.$
Let $\I = \{I\subseteq [m+n]|1, m+n\in I\}$. We consider $I\in\I$ as 
a strictly descending sequence of integers as in Equation \eqref{dod}, and define the corresponding element $f_{I} \in U(\fn^-)$ as in \eqref{ddd}.
We use the set $\I$ to index an expression for the \v Sapovalov element
$\gth_\eta \in U(\fb^-)$ analogous to that given in Theorem \ref{bb}.
For $k\in [m+n-1],$ consider the set
\[\cS_k=\{I\subseteq \I|k,k+1 \in I   \}.\]
\noi
Define as before,
\be r(I) = \{s-1|s \in \bar I\} \ee
where $\bar I$ be the complement of $I$ in $\I$.
Also let  $I^\pm$ be as in \eqref{qrt} 
If $k=m+n-1$, then $I^- \notin \I$ and if $k=1,$ then $I^+ \notin \I$. As before, set
$f_{I^-}=0$ or $f_{I^+}=0$ in these cases. Next set
\be \label{2tam}
 h_i=\left\{ \begin{array}
  {cc}  h_{\gs_i }  + (\gr,{\gs_i })-1  &\mbox{if} \;\; i \in [m-1]
\\  h_{\gt_{i+1-m }} + (\gr,{\gt_{i+1-m }}) &\mbox{if} \;\; m\le i \le m+n-2.
\end{array} \right. \ee
Then we have an analog of Lemma \ref{flo} (b),
\be \label{lan} h_{m+j-1} -h_{m+j}
 = h_{\gc_{j}}-1 \mbox{ for } j\in [n-2].\ee
If we  set  $h_{m+n-1}=0$, then \eqref{lan} holds also for $j= n-1$. 
Now define $$H_J=\prod_{i \in r(J)}h_{i}.$$

\bt\label{bsb1} Suppose  $\eta =\gep_1 -\gd_{n}$, $v$ is a highest weight vector in the $\fg$-Verma module $M(\gl)$ of weight $\gl$, and set
\[\Gt_\eta= \sum_{J\subseteq \I} f_J H_J.\]
Then $\Gt_\eta v$ is a $\fg_0$-highest weight vector. If $(\gl+\gr,\eta)=0,$  then $\Gt_\eta v$ is a $\fg$-highest weight vector, and so
$\Gt_\eta$ is a \v Sapovalov element for the pair $(\eta,1)$.
\et

\bpf This follows the same general pattern as the proof of Theorem \ref{bb}, but there are several differences resulting from Remark \ref{5r}, so we give some details. Suppose $I\in \cS_k$.
First if $\ga = {\ga_k}$ for $k\in [m-1]$ the proof of \eqref{aa}
shows that $$e_\ga (f_{I}H_I + f_{I^+}H_{I^+} + f_{I^-}H_{I^-})v = 0.$$
Next suppose that $k=m+j$ for $j \in [n-1]$, and set $\ga=\gc_j$, see
\eqref{heb}. Instead of \eqref{ay}, we set
\be \label{az}I_1 = \{m+n, i_1, \ldots i_{g-1},k+1 \}\mbox{ and } I_2 = \{k, i_{g+2}, \ldots i_s, 1 \}. \ee Note that $e_{-\ga}$ is a factor of $f_I$. Thus from Remark  \ref{5r}, we have
instead of \eqref{2vam},
\be \label{ata}
 e_\ga f_{I}v = -f_{I_1} h_\ga f_{I_2}v = f_{I_1} f_{I_2}(1-h_\ga)v. \ee
By  \eqref{dog22}, we have
$H_{I^+} = h_{m+j-1}H_I$ and $H_{I^-} = h_{m+j}H_I$.
Note that Lemma \ref{hog} still holds in this situation. Now in place of \eqref{aa} we have using \eqref{ata},
\be \label{asa} e_\ga (f_{I}H_I + f_{I^+}H_{I^+}+ f_{I^-}H_{I^-})v = f_{I_1}f_{I_2} (1-h_{\gc_j} +h_{m+j-1} -h_{m+j} ) H_Iv\ee
and this is zero by \eqref{lan}. This proves the first statement of the Theorem.\\ \\
Finally suppose that $k=m$ so that $\ga$  is the odd simple root and
$e_\ga = e_{m,m+1}$. Fix $I \in \cS_m$. Then
$f_{I^+}$ (resp. $f_{I^-}$) contains the odd root vector $e_{m+1, i_{g+2}}$ (resp.
$e_{ i_{g-1},{m}}$) as a factor. Hence 
$$e_\ga f_{I^+}H_{I^+}v= f_{I_1}f_{I_2}h_{m-1} H_Iv  \mbox{ and } e_\ga f_{I^-}H_{I^-}v= f_{I_1}f_{I_2}h_{m} H_Iv.$$ On the other hand 
$$e_\ga f_{I}H_{I}v= f_{I_1}h_\ga f_{I_2} H_Iv = f_{I_1}f_{I_2}(1+h_{\ga}) H_Iv.$$
Therefore by \eqref{2tam}, we have for $\gl\in \cH_\eta$, 
\by \label{as9} e_\ga (f_{I}H_I + f_{I^+}H_{I^+} + f_{I^-}H_{I^-})v &=& f_{I_1}f_{I_2} (1+h_{\ga} +h_{m-1} +h_{m} ) H_Iv\nn\\ &=&
 f_{I_1} f_{I_2}(h_\ga+ h_{\gs_{m-1}} +h_{\gt_{1}} + (\gr,{\gs_{m-1}} +{\gt_{1}}))H_Iv.\nn\\
&=& f_{I_1} f_{I_2}(h_\eta+  (\gr,\eta))H_Iv.\nn\\
&=& 0. \nn
\ey
\epf


\section{\v Sapovalov elements for arbitrary Borel subalgebras.} \label{sab}
Let $\fb^\dist$  denote the distinguished Borel subalgebra of $\fgl(m,n)$.  
We denote the set of Borel subalgebras with the same even part as $\fb^\dist$ by $\cB$.
For $\fb\in\cB$ and $\eta$ a positive root of both $\fb^\dist$ and $\fb$ we describe the \v Sapovalov elements with respest to $\fb$.
In \cite{M} 3.3, the set $\cB$ is described in terms of shuffles.  Here the notation is slightly different. We write a permutation $\gs$ of the set
 $\|m+n\|:=\{1,2, \ldots,m, 1',2',\ldots,n'\}$ in one-line notation as
\[ \underline{\gs} = (\gs(1), \gs(2), \ldots ,\gs(m), \gs(1'),\gs(2'),\ldots,\gs(n')).\] The last $n$ entries of $\|m+n\|$ are called {\it primed}, the others {\it unprimed}.  The shuffle condition is that both $
1,2, \ldots,m$ and $ 1',2',\ldots,n'$ are subsequences of $\underline{\gs}$. 
A useful way to think  about this is by using the Dynkin-Kac diagram of $\fb$, see \cite{M} section 3.4. The simple roots $\ga_{1}, \ldots, \ga_{m+n-1}$
are used to label the nodes from left to right.  Then augment the diagram by adding  the entries of 
$\underline{\gs}$, in order  above the diagram between the nodes, and at both ends. Without loss we assume that $\eta = \gep_1-\gd_n$ is the longest odd positive root of $\fg=\fgl(m,n)$, since any odd root which is positive in both 
$\fb^\dist$  and  $\fb$ is the longest root in a  possibly smaller Type A superalgebra.  This assumption implies that   
$\gs(1)= 1$ and $\gs(n')=n'$. 
From the diagram we can immediately read off the simple root corresponding to each node.
Each simple root $\ga_i$ lies between a pair of neighbors in $\underline{\gs}$, and there are 4 posibillities. The rule is
\by \label{ab3}\mbox{ if the neighbors of } \ga_i \mbox{ are }&& (a, a+1),\;\; (a, b'),\;\;\; (a',b),\;\;\; (b', b+1'),\nn\\
\mbox{  then the root } \ga_i \mbox{ is respectively }&& 
 {\gep_a}- {\gep_{a+1}},\;\;  {\gep_{a}} - {\gd_{b}},\;\;   {\gd_{a}} -
{\gep_{b}},\;\;{\gd_{b}}-{\gd_{b+1}}.  \nn\ey

\bexa 
{\rm We give an example of an augmented  Dynkin-Kac diagram when $m=4, n=3$. } 

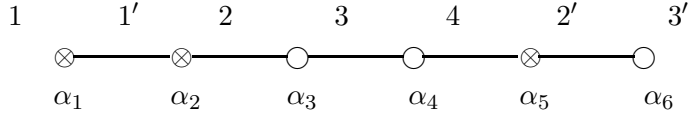
\vspace{0.5cm}

\begin{picture}(60,40)(-73,-20)
\thinlines 
\put(-6.7,-0.3){$\otimes$}
\put(-6.7,-15){$\ga_1$}
\put(1,3){\line(1,0){36.2}}
\put(46,3){\line(1,0){35.8}}
\put(37.7,-0.3){$\otimes$}
\put(37.7,-15){$\ga_2$}
\put(-23.7,15){$1$}
\put(17.7,15){$1'$}
\put(55.7,15){$2$}
\put(99.7,15){$3$}
\put(141.7,15){$4$}
\put(183.7,15){$2'$}
\put(225.7,15){$3'$}
\put(46,3){\line(1,0){35.8}}
\put(85.7,2.3){\circle{8}}
\put(81.7,-15){$\ga_3$}
\put(90,3){\line(1,0){36}}
\put(129.7,2.3){\circle{8}}
\put(127.7,-15){$\ga_4$}
\put(133,3){\line(1,0){36}}
\put(216.7,-15){$\ga_6$}
\put(169.7,-0.3){$\otimes$}
\put(169.7,-15){$\ga_5$}
\put(177,3){\line(1,0){36}}
\put(216.7,2.3){\circle{8}}
\end{picture}
\eexa  \noi
In this example $\underline{\gs} = (1, 1', 2,3,4,2',3')$ and 
$$\ga_1 =\gep_1- \gd_ 1,\;\ga_2= \gd_1-\gep_2,\;\ga_3= \gep_2-\gep_3,\;\ga_4= \gep_3-\gep_4,\;\ga_5= \gep_4-\gd_2,\;\ga_6=\gd_2-\gd_3.$$
In a similar way we can read off any positive root from the diagram, and it follows that $\sum_{k=1}^{m+n-1} \ga_k= \eta$.
Under the isomorphism $\fh^*\lra \fh$, suppose $\ga_k$ maps to $h_k$. We deduce that 
\be  \label{cy6} \sum_{k=1}^{m+n-1} h_k= h_\eta.\ee
Let $\I = \{ I \subseteq \|m+n\| | 1, n'\in I\}$. 
The shuffle $\gs$ induces a decreasing order on the set $\|m+n\|$. Explicitly
\[ \gs(n')=n' >  \gs(n-1') > \ldots
> \gs(2) > \gs(1)=1.\] 
and for $I\in \I$, $J \subseteq I$, we use this  order to form the products $f_J$ by analogy with Equation \eqref{ddd}.  
\\ \\
Now fix a simple root $\ga=\ga_k$ and  let $i(k)$ (resp. $j(k)$) be the number of grey nodes to the left (resp. right) of $\ga$ in the diagram.  
\noi  If $a$ and $b$ are the left and right neighbors of $\ga$ respectively, define $$S_\ga = \{I\in \I| a, b \in I\},$$ and 
note that $a$ and $b$ can be primed or unprimed numbers. 
 For $I\in \I$ we can write 
\be \label{az1} I =\{n', \ldots c,b, a,*, \ldots , 1 \}\ee 
 as a disjoint union, $ I_1\cup I_2$, where
\be \label{ay1}I_1 = \{n', \ldots ,b\}\mbox{ and } I_2 = \{a,\ldots , 1 \}.\ee
The bilnear form $(\;,\;)_{m,n}$ on $\fh^*$ defined in \eqref{hcb} induces an isomorphism  $\fh^* \lra \fh$ with $\ga \mapsto h_\ga$.
On the other hand there is  a bilinear form on $\fh$ induced by the supertrace Str, that is $(h,h')\mapsto$ Str$(hh')$. These are related by 
\be \label{3vab} (\ga, \gb) = \mbox{ Str}(h_\ga h_\gb). \ee
We need an explicit expression for $h_k = h_{\ga_k},$ in terms of the neighbors of $\ga_k =\ga$.
For $\ga$  is even, we already know the answer, 
\be \label{3vaa} h_{k} =(-1)^{i(k)}(e_{aa} - e_{bb}).\ee
Indeed $\ga$  is a root of $\fgl(m)\op 0$ if $i(k)$ is even, and a root of $0\op \fgl(n)$ if $i(k)$ is odd.
Thus \eqref{3vaa} follows from Remark \ref{5r}.  

\bl \label{yay1} We have
\bi \itema 
\be \label{3vza} [e_\ga,e_{- \ga} ]=(-1)^{i(k)}h_k.\ee
\itemb 
If $\ga$  is odd, then 
\be \label{3vaa1} h_{k} =(-1)^{i(k)}(e_{aa} + e_{bb}).\ee
\ei \el \bpf 
(a) This is an easy computation. (b) It seems that the easiest way to proceed is to take \eqref{3vaa1} as the definition, and then show that \eqref{3vab} holds.  
Consider  the Cartan matrix $(\ga_{k}, \ga_j)$. Suppose first that $k \neq 1, m+n-1$. The row indexed by $\ga_k$
has consecutive entries $-1, 0, 1$ (resp. $1, 0, -1$) if  $i(k)$ is even (resp. odd) with 0 on the main diagonal. All other entries in this row are zero. If $k=1$ (resp. $k=m+n-1$), the first (resp. last) entry in the above sequences must be deleted, since  0 is on the main diagonal.
 We show that the same holds for the Gram matrix of the  bilinear form on $\fh$ induced by the supertrace form, $(h,h')\mapsto$ Str$(hh')$. For example if $ {i(k)}$ is even,  then $a$ is unprimed. If $\ga_{k-1}$ is even,  then $h_{k-1} = e_{**} - e_{aa} $ where 
 $*$ is the entry in $I$ to the right of $a$, see \eqref{az1}.  
Then Str$(h_{k-1}h_{k})  =-1.$ 
If $\ga_{k-1}$ is odd,  then $h_{k-1} = 
(-1)^{i(k-1)}
(e_{**} + e_{aa}) $ and $i(k-1)$ is odd. Again we obtain  Str$(h_{k-1}h_{k})  =-1.$ Similarly Str$(h_{k}h_{k+1})  =1.$ We leave the cases  where 
$i(k)$ is odd,  $k=1$  or  $k=m+n-1$ to the reader.
\epf\noi  Note also that 
\be \label{4vaa} 
(e_{aa} \pm e_{bb})f_{I_2} =f_{I_2}(e_{aa} \pm e_{bb}+1) \mbox{ if }  I_2  \mbox{ is not a singleton}. \ee
\bl \label{goog} Suppose that $v$ is a highest weight vector, and $I\in S_\ga$. Then
\bi \itema
If $\ga$ is even we have
\by \label{2vbm}
 e_\ga f_{I}v &=& f_{I_1}[e_\ga,e_{- \ga} ]f_{I_2}v \nn \\ &=&
\left\{ \begin{array}
  {cc} (-1)^{i(k)}f_{I_1}h_kv&\mbox{if} \;\; I_2 = \{1\} \mbox{ is a singleton }
\\ f_{I_1}f_{I_2} (1+(-1)^{i(k)}h_k)v &\mbox{ otherwise}.
\end{array} \right. \ey

\itemb
If $\ga$ is odd we have
\by \label{2vcm}
 e_\ga f_{I}v &=&(-1)^{j(k)} f_{I_1}[e_\ga,e_{- \ga} ]f_{I_2}v \nn\\ &=&
\left\{ \begin{array}
  {cc} f_{I_1}h_kv&\mbox{if} \;\; I_2 = \{1\} \mbox{ is a singleton }
\\ f_{I_1}f_{I_2} (h_k+(-1)^{i(k)})v &\mbox{ otherwise}. 
\end{array} \right. \ey

\ei
\el 
\bpf 
(a) The first equality follows since 
$e_\ga$ commutes with  $f_{I_1}$ and then the second follows from \eqref{3vaa} and \eqref{4vaa}.\\ 
\\
(b) If $\ga$ is odd then $e_\ga$ anticommutes with $j(k)$ factors in $ f_{I_1}$ and commutes with the others. Note that the parity of  ${j(k)}$ depends only on the first and last entries in $I_1$.  This gives the first equality. 
The condition that   
$\gs(1)= 1$ and $\gs(n')=n'$  implies that   
the total number of odd simple roots is odd, so if $\ga$ is odd, then 
$i(k)+j(k)$ is even. Hence the second equality holds by Lemma \ref{yay1} (a) and \eqref{4vaa}.

\epf \noi 
Next for $I\in S_\ga$, set
\be\label{qot} I^+ = I \backslash \{a\}, \mbox{ and } I^- = I \backslash \{b\}.\ee 
If $I^\pm \notin \I$, then as before, we set 
$f_{I^\pm}=0$.
Then we have an analog of Lemma \ref{hog}.
\bl \label{hoog} We have \bi \itema 
If $\ga$ is even, then  $e_\ga f_{I^+}v = f_{I_1} f_{I_2}v$ and $e_\ga f_{I^-}v = -f_{I_1} f_{I_2}v$. 
\itemb If $\ga$ is odd, then
\by \label{3vw}
 e_\ga f_{I^-}v =-
(-1)^{i(k+1)}f_{I_1}f_{ I_2}v,\ey
and 
\be \label{3vbm}
 e_\ga f_{I^+}v = (-1)^{j(k)}f_{I_1}f_{I_2}v  = (-1)^{i(k)}f_{I_1}f_{I_2}v.  \ee
\ei  
Furthermore Lemma \ref{hog} $ ($c$)$ holds.
\el
\bpf (a) If $\ga$ is even, then ``even rules" apply to the  commutator of 
$e_\ga$ and any root vector which is a factor of $ f_{I^\pm}$.  Thus we get the same result as Lemma \ref{hog} (a).\\
(b) We prove \eqref{3vw}.
Define 
\be \label{ay3}I^-_1 = \{n', \ldots ,c\}\mbox{ and } I^-_2 = \{a,\ldots , 1 \} =I_2 \ee
Then $f_{I^-} = f_{I^-_1} e_{ca} f_{I^-_2}$, so 

\by \label{3Aw}
 e_\ga f_{I^-}v &=&   (-1)^{j(k+1)}f_{I^-_1} e_{ab} e_{ca} f_{I^-_2}v\nn\\
&=&\left\{ \begin{array} {cc} (-1)^{j(k+1)}f_{I_1}f_{I_2}v&  \mbox{ if }  \ga_{k+1} \mbox{ is even  }
\\ (-1)^{j(k+1)+1}f_{I_1}f_{I_2}v &\mbox{ if } \ga_{k+1} \mbox{ is odd}.
\end{array} \right. \nn \ey
The first equality holds since there are exactly ${j(k+1)}$ factors in $f_{I^-_1}$ that are odd root vectors. 
For the second, note that if $  \ga_{k+1} \mbox{ is even}$ (resp. odd), then 
$e_{ca}$ is an odd (resp. even) root vector, and hence $[e_{ab}, e_{ca}] = e_{cb}$ (resp.  $[e_{ab}, e_{ca}] = -e_{cb}$).
 Now { if }  $\ga_{k+1} \mbox{ is even}$, then mod 2, 
${j(k+1)} \equiv 
{i(k+1)}+1$, while { if } $\ga_{k+1} \mbox{ is odd}$, then ${j(k+1)}+1 \equiv 
{i(k+1)}+1$. This gives \eqref{3vw}.
\\ \\
 The proof of  \eqref{3vbm} 
is similar to \eqref{3vw} but easier.  Note that, if 
\be \label{ay2}I^+_1 = \{n', \ldots ,b\}=I_1 \mbox{ and } I^+_2 = \{*,\ldots , 1 \},  \ee
we have $f_{I^+}v = f_{I^+_1} e_{b*} f_{I^+_2}v$.
The last statement in the Lemma is easy to prove.\epf \noi Recall $h_k$ defined by \eqref{3vaa} or \eqref{3vaa1}, and
 set  $s_{k} =h_1 +\ldots+ h_k= s_{k-1} + h_k$. Next let $$\ell(k) = | i \in \{[ k-1]| \ga_i \mbox{ is a root of 
} 0\op \fgl(n)
   \}|- | i \in \{[k-1]| \ga_i \mbox{ is a root of }  
\fgl(m)\op 0
 \}|,$$
and define $\bar i(k) = (1-(-1)^{i(k)})/2.$ Thus $\bar i(k) \in \{0, 1\}$, and $\bar i(k) \equiv i(k)$ mod 2. Set $d_k = \ell(k) +\bar i(k)$. 
 Thus $d_1=0 $.
\bl  \label{by2}
 If  $k\in [m +n-2]$,  we have the following recurrence,  
$$d_{k+1} = d_{k} -(-1)^{i(k+1)}.$$
\el
\bpf 
We have to show that
\be  \label{cy2} \ell(k+1) +\bar i(k+1)= \ell(k)+\bar i(k) -(-1)^{i(k+1)}. \ee
If $\ga_k$ is even, then $i(k+1)= i(k).$ Thus \eqref{cy2} becomes
\be  \label{cy3} \ell(k+1) = \ell(k) -(-1)^{i(k)}. \ee
Now if $i(k)$ is even (resp. odd), then $\ga_k$  is a root of 
$\fgl(m)\op 0$ 
(resp. $0\op \fgl(n)$) so $ \ell(k+1) = \ell(k) -1$  (resp. $ \ell(k+1) = \ell(k) +1$). This gives \eqref{cy3}. \\ \\
 If $\ga_k$ is odd, then 
$\ell(k+1) = \ell(k),$ and $i(k+1)= i(k)+1,$ so we have to show
\be  \label{cy4} \bar i(k+1)= \bar i(k) +(-1)^{i(k)}.\nn\ee
This follows easily by considering the cases, $i(k)
$ even and odd separately.
\epf
\bl  \label{cy7} $d_{m+n-1} =-(\gr,\eta).$
\el 
\bpf Let  $r$ (resp. $s$) be the total number of roots  in the Dynkin-Kac diagram belonging to $0\op \fgl(n)$
  (resp. 
$\fgl(m)\op 0$). Thus $(\gr,\eta) = -r+s$ by Remark \ref{5r} and  Lemma \ref{flo} (a).  
If $\ga= \ga_{m+n-1}$ is even, then it is a root of $0\op \fgl(n)$. Thus 
 $\ell({m+n-1}) = r-1-s$. Also  $i({m+n-1})$ is odd, so $
d_{m+n-1} = \ell({m+n-1}) +1 = r-s$. 
If $\ga$ is odd, then $i({m+n-1})$ is even, so $d_{m+n-1} = \ell({m+n-1}) = r-s$. This gives the result.
\epf  \noi 
If $(a,b)$ are the neighbors of $\ga_k$, and $1 < k<m+n-1$,
set

 \be  \label{cy8} h_{\hat a} =(-1)^{i(k)}(s_{k-1} -d_{k-1}) \mbox{ 
and } h_{\hat b}  = (-1)^{i(k+1)}(s_{k}-d_{k}).\ee 
If $k= 1$ (resp. $k= m+n-1$), then $I^+\notin \I$ (resp. $I^-\notin \I$), but we still define $h_{\hat b}$ (resp. $h_{\hat a}$) as in \eqref{cy8}. 
Note that $b$ is also the left neighbor of $\ga_{k+1}$, and Equation
\eqref{cy8} is  consistent with this fact. 
Also  \eqref{cy8} is very convenient for the proof of the main result because of \eqref{cz1} below. However we need another ingredient to define the 
coefficients $H_I$. If $e$ is the right  neighbor of $\ga_j$, then set $\ttr(e) =j$. Then define $$t_e = (-1)^{i(\ttr(e)+1)}(s_{\ttr(e)}-d_{\ttr(e)}),$$ and $H_I= \prod_{e \in \bar I} t_e$, where $\bar I$ is the complement of $I$ in $\|m+n\|$. 
For example  $\ttr(a) = k-1, \ttr(b) = k$, so by  
\eqref{cy8} $t_a = h_{\hat a}$ and  $t_b = h_{\hat b}$.   
If $(a,b)$ are as above, it follows from \eqref{qot} that 
 \be \label{cz1}
H_{I^+} =   h_{\hat a} H_{I} \mbox{ and } H_{I^-} =h_{\hat b} H_{I}.\ee


\bt\label{bsb9} Set
\[\Gt_\eta= \sum_{J\subseteq \I} f_J H_J.\] Then 
$\Gt_\eta$ is a \v Sapovalov element $\gth_\eta$ for the pair $(\eta,1)$ for the Borel subalgebra $\fb$.
\et
\bpf  
Assume first that   $k\neq 1, m+n-1.$ 
 If $\ga$ is even then, 
 ${i(k+1)}= {i(k)}$. Hence 
by \eqref{2vbm}, \eqref{cz1} and Lemma \ref{hoog} (a), then \eqref{cy8} and Lemma \ref{by2},
  \by \label{ap} e_\ga( f_{I}H_I + f_{I^+}H_{I^+} + f_{I^-}H_{I^-})v &=& f_{I_1}f_{I_2} (1+(-1)^{i(k)}h_{k}+h_{\hat a} -h_{\hat b} ) H_I v\nn\\
&=&(-1)^{i(k)} f_{I_1}f_{I_2} (h_{k}+s_{k-1} -s_{k}) H_I v\nn\\
&+& f_{I_1}f_{I_2} (1-(-1)^{i(k)}d_{k-1}  +(-1)^{i(k)}d_{k} ) H_I v\nn\\
&=&0.\ey
If $k=1$, then $I_2$ is a singleton and $d_{1}  =i(2) =0,$ so by \eqref{2vbm}, 
\eqref{3vbm}  and  \eqref{cy8}, 
\by \label{ap7} e_\ga( f_{I}H_I + f_{I^-}H_{I^-})v &=& f_{I_1} (h_{1}-h_{\hat b} ) H_I v\nn\\
&=&f_{I_1}(h_{1}-s_{1}+d_{1} ) H_I v\nn\\
&=&0.\nn\ey
Now assume  $\ga$ is odd and $k\neq 1, m+n-1.$ 
 Then, by Lemma \ref{goog} (b),  \eqref{3vbm} and  
 \eqref{3vw}, \eqref{cz1}, then 
\eqref{cy8} and
Lemma \ref{by2},
\by \label{aq}e_\ga( f_{I}H_I + f_{I^+}H_{I^+} + f_{I^-}H_{I^-})v &=& f_{I_1}f_{I_2} (h_{k}+(-1)^{i(k)} +(-1)^{i(k)}h_{\hat a} -(-1)^{i(k+1)}h_{\hat b} ) H_I v\nn\\
&=&
 f_{I_1}f_{I_2} ((h_{k}+s_{k-1} -s_{k}) H_I v\nn\\
&+& f_{I_1}f_{I_2} ((-1)^{i(k)} -d_{k-1} +d_{k} )
H_I v\nn\\&=&0.\ey
If $k=1,$ then by \eqref{2vcm}, \eqref{3vw},   \eqref{cy8} and \eqref{cz1}, 
\by \label{a1q}e_\ga( f_{I}H_I + f_{I^-}H_{I^-})v &=& f_{I_1} (h_{1}-(-1)^{i(k+1)}h_{\hat b} ) H_I v\nn\\
&=&
 f_{I_1} ((h_{1}-(s_{1}- d_{1}))H_I v\nn\\&=&0.\nn\ey
Now suppose that $k=m+n-1$, and that $v$ has weight $\gl\in \cH_\eta$. 
\noi
 If $\ga$ is even, then $i(k)$ is odd. In the sequence of equalities below the first follows from \eqref{2vbm} and \eqref{cz1}, the second from \eqref{cy8}. 
Then the  third follows from \eqref{cy6} and Lemma  \ref{by2}, the fourth from Lemma \ref{cy7} and the final equality holds since $\gl\in \cH_\eta$. 
\by \label{ap9} e_\ga( f_{I}H_I + f_{I^+}H_{I^+})v 
&=& f_{I_2} (1-h_{k}+h_{\hat a}) H_I v\nn\\
&=&f_{I_2} (1-h_{k}-s_{k-1} +d_{k-1}) H_I v\nn\\
&=& - f_{I_2} (h_{\eta} - d_{m+n-1}) H_I v \nn\\
&=&- f_{I_2} (h_{\eta}+(\gr,\eta)) H_I v= 0.\nn\ey
 Suppose $\ga$ is odd. Then  $i(k)$ is even. For the first two  equalities 
below  we use \eqref{2vcm},  
\eqref{3vbm},  \eqref{cz1} and \eqref{cy8}. 
The remainder of the proof is similar to the case where $\ga$ is even.
\by \label{aq5}e_\ga( f_{I}H_I + f_{I^+}H_{I^+})v
 &=& f_{I_2} ((-1)^{i(k)}+h_{k}+(-1)^{i(k)}h_{\hat a}) H_I.\nn\\
 &=& f_{I_2} (1+h_{k}+s_{m+n-2}-d_{m+n-2}) H_I \nn\\
 &=& f_{I_2} (h_{\eta} - d_{m+n-1}) H_I \nn\\
&=& f_{I_2} (h_{\eta}+(\gr,\eta)) H_I v= 0.\nn\ey
\epf


\section{\v Sapovalov  elements as determinants of Hessenberg matrices.}\label{shm}
\subsection{Hessenberg  Matrices.} 
An $n\ti n$ matrix $B=(b_{ij})$ is {\it  $($upper$)$ Hessenberg  of order $n$}
if $b_{ij}=0$ unless $i\le j+1$.
The only assumption on the Hessenberg  matrix $B$ is that the entries on the subdiagonal commute with all other entries.
We define a noncommutative determinant of the $n\ti n$ matrix $B=(b_{ij})$,
working from left to right, by
\be \label{lr}{\stackrel{\longrightarrow }{{\rm det}}}(B) =  \sum_{w \in \mathcal{S}_n} sign(w) b_{w(1),1} \ldots b_{w(n),n}, \ee
Cofactor expansions of ${\stackrel{\longrightarrow }{{\rm det}}}(B) $ are valid as long as the overall order of the terms is unchanged.
\bl \label{hes}
Suppose that $B$ is Hessenberg of order $n$.
\bi \itema For a fixed $q\in [n-1],$ let $\ttT = -b_{q+1q}.$ Then
\be \label{sam2}{\stackrel{\longrightarrow}{{\rm det}}} B =
\ttT  {\stackrel{\longrightarrow}{{\rm det}}}B'' + {\stackrel{\longrightarrow}{{\rm det}}}B',\ee
where  $B'$ is obtained   from  $B$ by setting $\ttT =0,$
and  $B''$ is obtained   from  $B$ by deleting the row and column containing $\ttT$.
\itemb The matrix $B'$ is block upper triangular, with two diagonal
blocks which are upper Hessenberg of order $q$ and $n-q$.
\itemc The matrix
$B''$ is upper Hessenberg of order $n-1$. Also any term in  the expression \eqref{lr} for ${\stackrel{\longrightarrow }{{\rm det}}}(B)$ which contains a factor  of the form $b_{iq}$ or $b_{q+1j}$ 
 cannot occur in ${\stackrel{\longrightarrow}{{\rm det}}} B''$.
\end{itemize}
\el
\bpf Part (a) follows by separating the products in \eqref{lr} that contain $\ttT$ from those that do not. Note that $\ttT$ commutes with all entries in $B$, and that the order of all other factors of the products is unchanged.  The rest is easy. \epf \noi
We consider determinants of Hessenberg matrices with entries on or above the diagonal  from $\fn^-$ and subdiagonal entries which commute with all other matrix entries. Explicit expressions in $U(\fb^-)$ for \v Sapovalov elements   can
be obtained as complete expansions of suitable determinants of this kind.
There are significant differences in the complete expansions
 depending on
the ordering of entries in the matrices.
We consider two orderings on the set of positive roots, and explain
the relationships between the  determinants without reference to \v Sapovalov elements.  First consider the matrix with entries in $U(\fgl(m))$, or any $U(\fgl(r,s))$ with $r+s=m$.

\be\label{ccc1}
D=  \left[ {\begin{array}{ccccc}
e_{m,m-1}&e_{m,m-2} & \hdots &e_{m ,2}  & e_{m,1} \\
 -a_{m-2} &e_{m -1,m-2} & \hdots&e_{m -1,2} & e_{m -1,1} \\
0& -a_{m-3} & \hdots &   e_{m -2,2} &e_{m -2,1} \\
 \vdots  & \vdots & \ddots & \vdots& \vdots \\
0&0&  \hdots &-a_{1} & e_{2,1} 
\end{array}}
 \right].
\ee
\noi
Note that in \eqref{ccc1} the entries from the same row from $\fn^-$ have the form $e_{i,*}$ for some fixed $i$.
Another possibility is to require that all entries from the same row have the form $e_{*,i,}$. 
 Thus consider the matrix, where $c_i = a_i+1$,


\be\label{coc}
E=  \left[ {\begin{array}{ccccc}
e_{2,1}&e_{3,1} & \hdots  & e_{m-1,1}  & e_{m,1} \\
 -c_{1} &e_{3,2} & \hdots & e_{m-1,2}  & e_{m ,2} \\
0& -c_{2} & \hdots& e_{m -1,3}&    e_{m ,3} \\
 \vdots  & \vdots & \ddots & \vdots   & \vdots \\
0&0& \hdots & -c_{m-2} & e_{m,m-1} 
\end{array}}
 \right].
\ee
The next result is used in the proof of Theorem \ref{stc}.
\bp \label{2aa}We have
\[ {\stackrel{\longrightarrow }{{\rm det}}}\;D = {\stackrel{\longrightarrow }{{\rm det}}}\;E\]
\ep
\bpf This is proved using  cofactor expansion and induction. Rather than giving full details, an  example might be more helpful here. It is easy to see how the example generalizes. In the  general case, we  need a commutation relation for cofactors of $E$.\epf
\bl \label{ee} Let $E^1$, $E^2$ be the cofactors of entry $e_{m,m-1}$ and $-c_{m-2}$ in $E$ respectively. Then $[{\stackrel{\longrightarrow }{{\rm det}}}\;E^1,e_{m,m-1}]= -{\stackrel{\longrightarrow }{{\rm det}}}\;E^2$.
\el
\bpf  Note that
$E^1$ is obtained from $E^2$ by replacing the entries $e_{m,i}$ with $e_{m-1,i}$ 
for $i\in [m-2]$.  Since $i<m-1$, it is impossible for both $e_{m-1,i}$ and  $e_{m,m-1}$  to be odd.  Thus $[e_{m-1,i}, e_{m,m-1}] = -e_{m,i}$, the result follows.
\epf

\bexa
{\rm Let
\be
E=  \left[ {\begin{array}{ccc}
e_{2,1}&e_{3,1} &  e_{4,1} \\
 -c_1 &e_{3,2} & e_{4,2} \\
0& -c_{2} &     e_{4,3}
\end{array}}
 \right].\nn
\ee
By cofactor expansion along the last row, 
\by {\stackrel{\longrightarrow }{{\rm det}}}\;E
&=& {\stackrel{\longrightarrow }{{\rm det}}}\;
E^1 e_{4,3}
 + c_2{\stackrel{\longrightarrow }{{\rm det}}}\;
E^2.\nn
\ey
where
\[E^1 = \;\left[ {\begin{array}{ccc}
e_{2,1}&e_{3,1}  \\
 -c_1 &e_{3,2}
\end{array}}
 \right] \mbox{ and  }
 E^2= \left[ {\begin{array}{ccc}
e_{2,1}&  e_{4,1} \\
 -c_1 & e_{4,2}
\end{array}}
 \right].
\]
\noi
Next set
\[D^1 = \;\left[ {\begin{array}{ccc} e_{3,2} &  e_{3,1} \\
 -a_1 &
e_{2,1}
\end{array}}
 \right] \mbox{ and  }
 D^2= \left[ {\begin{array}{ccc}
e_{4,2}&  e_{4,1} \\
 -a_1 & e_{2,1}
\end{array}}
 \right].\]
Now by Lemma \ref{ee},
\[[{\stackrel{\longrightarrow }{{\rm det}}}\;
E^1, e_{4,3}] = -{\stackrel{\longrightarrow }{{\rm det}}}\;
E^2, \]
so 
\by \label{1bb}{\stackrel{\longrightarrow }{{\rm det}}}\;E&=& e_{4,3}{\stackrel{\longrightarrow }{{\rm det}}}\;
E^1
 + (c_2-1){\stackrel{\longrightarrow }{{\rm det}}}\;
E^2\nn\\
&=& e_{4,3}{\stackrel{\longrightarrow }{{\rm det}}}\;
D^1
 + a_2{\stackrel{\longrightarrow }{{\rm det}}}\;
D^2
\ey
using induction for the second equality. This is the cofactor expansion of ${\stackrel{\longrightarrow }{{\rm det}}}\;D$  down the first column.
}\eexa


\subsection{ Application to 
\v Sapovalov elements in Type A.}
Expansions of
\v Sapovalov elements in Type A were already given in \cite{M2} Section 9, using determinants of a certain Hessenberg matrices. However it seems unlikely that this method will generalize, because the determinant of Hessenberg matrix of order $m$ has $2^{m-1}$ terms, and outside of Type A, the number of partitions of a positive root is rarely a power of 2.
Nevertheless the use of Hessenberg matrices gives more insight in the Type A
case.   In particular they can be used to give expressions for  \v Sapovalov elements using different orderings on the set of positive roots.\\ \\
Consider the Lie algebra $\fgl(m)$ with simple roots
$\ga_i = \gep_{i}-\gep_{i+1}$ for $i\in [m-1]$.
 Let  $\eta' = \gep_{1}-\gep_{m}$ and  $\mu \in \cH_{\eta'}$.  We show that $\gth_{\eta'}$ can be expressed as a determinant of a certain Hessenberg matrix.
Recall the elements $$h_i= h_{\gs_i} +(\gr,{\gs_i}) -1$$ from \eqref{rt4}, for $i\in [m-1]$, Then let $
\cD^m(\mu)$ be the matrix $D$ from \eqref{ccc1} with 
\be \label{dz1}a_i= h_i(\mu) = (\mu +\gr,{\gs_i}) -1.\ee
\bt \label{hea}The \v Sapovalov element for $\eta'$ satisfies
\[\gth_{\eta'}(\mu) = {\stackrel{\longrightarrow}{{\rm det}}}\;\cD^{m}(\mu)\]
for all $\mu\in \cH_{\eta'}$.
\et \noi
We give two proofs of Theorem \ref{hea}.  A comparison of the two approaches is necessary to prove Proposition \ref{11.5}. Define
$H_J$ from \eqref{1q}.
In the first we show that the complete expansion of
the determinant 
${\stackrel{\longrightarrow}{{\rm det}}}\;\cD^{m}(\mu)$
is the evaluation of 
the element $\Gt_{\eta'}$ from Theorem \ref{bb}
at $\mu\in \cH_{\eta'}$. In other words we show 
\be \label{da1}{\stackrel{\longrightarrow }{{\rm det}}}\;\cD^{m}(\mu)= \sum_{J\subseteq \I} f_J H_J(\mu).\ee
In the second we use induction on $m$, cofactor expansion and the fundamental Lemma \ref{1768} below. This Lemma which is \cite{M} Lemma 9.4.3, is the basis for the proof of the existence  of \v Sapovalov elements in the general case. Essentially the same proof is given in \cite{H2} Section 4.13. 
\\ \\
{\it Theorem \ref{hea}: First Proof.} \label{fpf}
Consider the matrix $\cD^{m}(\mu)$ from Theorem \ref{hea}.
We obtain the complete expansion of
${\stackrel{\longrightarrow}{{\rm det}}}\;\cD^{m}(\mu)$.
Let $\I$ be  as  in \eqref{120}.
 Each term in the complete expansion 
is obtained by choosing  a non-zero product of elements from each column, with each row occurring exactly once.   
The product of the 
{\it chosen elements} lying  above the subdiagonal has the form  $f_I$ for some $I\in \I$.
The proof of \eqref{da1} is completed by the following Lemma.
\hfill  $\Box$
\bl \label{ece} The product $a_I$ of subdiagonal terms accompanying $f_I$ is given by
\be \label{dd} a_I=\prod_{i\in r(I)} a_{i} =\prod_{i\in r(I)} h_{i}(\mu) = H_I(\mu).\ee
 \el \bpf 
From the form of the matrix $D$ in \eqref{ccc1} we see that 
$e_{j,*}$ is a factor of $f_I$ iff $j \in I$ iff $-a_{j-1}$ a factor of $a_I$.  Thus 
the product of subdiagonal terms accompanying $f_I$ must be
\be \label{dd1} \pm \prod_{i\in r(I)} a_{i} = \pm \prod_{i\in r(I)}  h_{i}(\mu).\ee
\noi
Next we make a remark about the (non-commutative) determinant of a Hessenberg matrix $B$ of order $n.$ For a general $n\ti n$ matrix the complete expansion of the determinant is indexed by permutations $w$ from the symmetric group of degree $n$.  In Equation \eqref{lr} the corresponding term is zero unless $w(i) \ge i-1$ for $2\le i\le n$.  Now it is easy to see the following.
\bl
Assume that in the expression for  $\stackrel{\longrightarrow }{{\rm det}}(B) $ given in \eqref{lr}, the term indexed by $w$ is nonzero.
 Let $C = |\{i| w(i) = i-1\}|$.  Then $sign(w) = |(-1)^{|C|}.$\el
\bpf Each $i\in C$ corresponds to an inversion $(i,i-1)$ in $w$. There are no other inversions.\epf
\noi This is fortunate for us because the set $C$ corresponds to the elements on the subdiagonal of $B$ and these all have the form $-a_j$ for some $j$. This means that \eqref{dd1} can be replaced with
 \eqref{dd} as in the statement of Lemma \ref{ece}.
\epf \noi
{\it Theorem \ref{hea}: Second Proof.} \label{spf} 
\noi  The \v Sapovalov element $\theta_{\eta}$ can be constructed inductively using the next Lemma.
\begin{lemma}\label{1768}
  Suppose that  $\alpha$ is a simple root and set $\mu = s_\alpha \cdot  \gl, \; \eta' = s_\alpha\eta.$
Assume that
\begin{itemize}
 \itema $p = (\mu + \rho, \alpha^\vee) \in \mathbb{N}
\backslash \{0\}$ and $q = (\eta, \alpha^\vee) \in \mathbb{N}
\backslash \{0\}$
\itemb $\mu\in \cH_{\eta'}$, and consequently $\gl\in \cH_{\eta}$.
\end{itemize}
Then the evaluation of the \v Sapovalov  elements
$\theta_{\eta'},  \theta_{\eta}$  at $\mu$ and $\gl$ satisfy
\be \label{epqm}  e^{p + q}_{- \alpha}\theta_{\eta'}(\mu) = \theta_{\eta}(\gl) e^p_{- \alpha}.\ee
\el
\bpf See the sources cited above. \epf \noi
Consider the Lie algebra $\fgl(m+1)$ with simple roots
$\ga_i = \gep_{i}-\gep_{i+1}$ for $i\in [m]$, and set $\gs_i = \gep_{1}-\gep_{i+1}$.
 Let $\ga=\ga_m$, $\eta' = \gep_{1}-\gep_{m}$ and $\eta = \gep_{1}-\gep_{m+1} = s_\ga \eta'$.
Suppose  $\ga, \mu, {\eta}$ satisfy the hypotheses of Lemma \ref{1768}. Then $\gl = s_\ga \cdot\mu \in \cH_{\eta}$. 
Thus since $(\gl + \rho, \alpha^\vee) =-p$, it follows that
\[1 = (\gl + \rho, \eta)= (\gl + \rho, \gs_{m-1} +\alpha) = -p+(\gl + \rho, \gs_{m-1}), \] which implies that $p=(\gl + \rho, \gs_{m-1}) -1:= a_{m-1}$.
If $ i\in [m-2]$, then $(\gs_{i},\ga)=0$, it follows that the $a_i$ as defined in  \eqref{dz1} satisfy
\be\label{gg7}
a_i=(\gl + \gr,\gs_{i} ) -1,\ee
 which means that $\cD^{m}(\mu)$ fits into the lower right corner of $\cD^{m+1}(\gl).$ 
Then  $\cD^{m+1}(\gl)$ is obtained from $\cD^{m}(\mu)$ by first adding a column on the left, whose only non-zero entry is $-a_{m-1}$ in the first row, and then adding a  row so that
$e_{m+1,i}$ is directly above $e_{m,i}$  for $i\in [m-1]$, and with $e_{m+1,m}$  as first entry.  
\\ \\
Next introduce 
two matrices  $\cD_1, \cD_2$ by

\be \label{giv} \cD_1=
\left[ {\begin{array}{clccccc}
e_{m+1,m}&\vline&0&\hdots & 0\\
\hline
0 &\vline & & & \\
0&\vline&    &&    \\
\vdots &\vline&  & \cD^{m}(\mu)&   \\
0&\vline&&&\\
\end{array}}
 \right],
 \cD_2=\left[ {\begin{array}{clccccc}
0&\vline&e_{m+1,m-1}&\hdots & e_{m+1,1}\\
\hline
-a_{m-1} &\vline & & & \\
0&\vline&    &&    \\
\vdots &\vline&  & \cD^{m}(\mu)&   \\
0&\vline&&&\\
\end{array}}
 \right].
\ee
Observe that $e_{m+1,m}$ commutes with all entries in $\cD^{m}(\mu)$ except for those in the first row, and that
\be \label{art} e_{m+1,m}^{p+1}e_{m,i} = (e_{m+1,m}e_{m,i}+pe_{m+1,i})e_{m+1,m}^p.\ee
This  implies that

\be \label{are} e_{m+1,m}^{p+1}{\stackrel{\longrightarrow}{{\rm det}}}\;\cD^{m}(\mu) = ({\stackrel{\longrightarrow}{{\rm det}}}\;\cD_1 + {\stackrel{\longrightarrow}{{\rm det}}}\;\cD_2)e_{m+1,m}^p.\ee
In \eqref{are}, $\cD_2$ arises from the second term in the factor  $(e_{m+1,m}e_{m,i}+pe_{m+1,i})$ from  Equation \eqref{art}, and similarly $\cD_1$ corresponds to the first term.
But ${\stackrel{\longrightarrow}{{\rm det}}}\;\cD_1 + {\stackrel{\longrightarrow}{{\rm det}}}\;\cD_2$ is just the cofactor expansion of ${\stackrel{\longrightarrow}{{\rm det}}}\;\cD^{m+1}(\gl)$ down the first column. 
By induction
\[\gth_{\eta'}(\mu) = {\stackrel{\longrightarrow}{{\rm det}}}\;\cD^{m}(\mu).\]
Now comparing \eqref{epqm} and \eqref{are}, and using the fact that $e_{-\ga}$ is not a zero divisor in $U(\fg)$ it follows that
\[\gth_\eta(\gl) = {\stackrel{\longrightarrow}{{\rm det}}}\;\cD^{m+1}(\gl).\]
\hfill  $\Box$
\\ \\
We record how \v Sapovalov elements can be expressed using non-commutative determinants  in the case of $\fg= \fgl(m,n)$. First define 
\be \label{5t} B_i = \left\{ \begin{array}
  {cc}  (\gl+\gr, \gs_i) -1\ &\mbox{ for } i\in [m-1],\\
	(\gl+\gr,\gt_{i+1-m }) & \mbox{ for } \;\; m\le i \le m+n-2.
\end{array} \right. \ee
With $h_i$  as in \eqref{2tam}, we have 
\be \label{1w} h_i(\gl) = B_i \mbox{ for all } i\in [m+n-2].\ee
Next consider the determinant

\be \label{gav} \A(\gl) = \left[ {\begin{array}{ccccccc}
e_{m+n,m+n-1}&e_{m+n,m+n-2} & \hdots & \hdots &\hdots &e_{m+n,2}  & e_{m+n,1} \\
 -B_{m+n-2} &e_{m+n-1,m+n-2} & \hdots & \hdots &\hdots &e_{m+n-1,2} &  e_{m+n-1,1} \\
0& -B_{m+n-3} & \hdots &\hdots &\hdots & e_{m+n-2,2} &   e_{m+n-2,1} \\
 \vdots  & \vdots & \ddots & \vdots& \vdots  &\hdots & \hdots \\
0 & 0 &-B_{m} &e_{m+1,m}&\hdots& \hdots& e_{m+1,1}\\
0 & 0 &&-B_{m-1} &e_{m,m-1}&\hdots& e_{m,1}\\
 \vdots  & \vdots & \ddots & \vdots& \vdots  &\hdots & \hdots \\
0&0&  \hdots &\hdots &\hdots & -B_1& e_{2,1} \\
\end{array}}
 \right],
\ee 
Recall the expresssion for the  \v Sapovalov element $\gth_\eta$ from Theorem \ref{bsb1}.
\bt\label{7t} We have
\[{\stackrel{\longrightarrow }{{\rm det}}}\;\A(\gl)
= \sum_{J\subseteq \I} f_J H_J(\gl).\]
Thus ${\stackrel{\longrightarrow }{{\rm det}}}\;\A(\gl)$
is the evaluation of the \v Sapovalov element $\gth_\eta$ at $\gl \in \cH_\eta$.
\et
\bpf This follows from \eqref{1w},  Theorem \ref{bsb1} and the argument of Theorem \ref{hea} showing that $f_J$ is correctly paired with $H_J$ in the above sum.\epf

\subsection{\v Sapovalov elements for arbitrary odd roots.}\label{ipu}
Assume $\gc$ is an  odd root. Although as we remarked in Section \ref{se1}, to compute  \v Sapovalov element $\gth_\gc$ we can assume  $\gc$ is  the highest odd root
 of $\fgl(m,n)$, we later apply our results when this is not the case. 
Thus it will be conveneint to have and explicit expressions for a \v Sapovalov elements for an arbitrary odd root.\\ \\
Let $t= m+1-r$. We need the \v Sapovalov element for $\gc=\gep_r - \gd_s$ in $U(\fgl(m,n))$.
Let
$\fg'$ be the subalgebra     $\fg=\fgl(m,n)$ 
with rows and columns indexed by the set 
$\{r, \ldots , {m+s}\}$.  Note that $\fg'$ and  $\fg$ both  share the same odd simple root vector $e_\gb= e_{m,m+1}$.
Also $\fg'\cong \fgl(t,s)$, and 
$\gc$ is the longest odd positive root for the  subalgebra  $\fg'$.
The \v Sapovalov element  $\gth_\gc$ is easily found using the matrix
\be \label{gbv} \A^{r,s}(\gl) = \left[ {\begin{array}{cccccc}
e_{m+s,m+s-1}&e_{m+s,m+s-2} &\hdots &\hdots &e_{m+s,r+1}  & e_{m+s,r} \\
 -A_{m+s-2} &e_{m+s-1,m+s-2} &\hdots &\hdots &e_{m+s-1,r+1} &  e_{m+s-1,r} \\
 \vdots  & \vdots & \ddots&\hdots &\hdots & \hdots \\
 \vdots  & \vdots & \ddots&\hdots &\hdots & \hdots \\
0& 0&\ddots  & -A_{r+1}& e_{r+2,r+1} &   e_{r+2,r} \\
0&0&\hdots &  \hdots  & -A_r& e_{r+1,r} \\
\end{array}}
 \right],
\ee
Note that all entries in the above matrix belong to the subalgebra $\fg'$ of $\fg$.
Set
\be \label{67t}\gs_i =\gep_{r} - \gep_{r+i}, 
\mbox{ and }\go_i =\gd_{i+1-t }-\gd_{s}.\ee
Then we define $A_i$ by
\be \label{6t} A_{i+r-1} = \left\{ \begin{array}
  {cc}(\gl+\gr,\gs_{i})-1 &\mbox{ for } i\in [t-1].\\
	(\gl+\gr,\go_{i})& \mbox{ for } \;\; t\le i \le t+s-2.
\end{array} \right. \ee
We have $(\gr, \ga)$ = 1 or $-1$, respectively if $\ga$ is a simple root of $\fgl(t)\op 0$ or $0\op \fgl(s)$ and   $(\gr,\gb) = 0$ for the unique simple odd root, see Remark \ref{5r}. These are the only properties of $\gr$ we need, so there is no need to introduce the analog of $\gr$ for $\fgl(t,s)$. 
\noi
Let
\be \label{2w} \I = \{I\subseteq \{r, r+1,\ldots, m+s\}|r, m+s\in I\}.\ee
 Then for$ J \in \I$, define $f_J, r(J)$ as usual.
Define $H_{i+r-1}\in U(\fh)$ by
\be \label{9t} H_{i+r-1} = \left\{ \begin{array}
  {cc}h_{\gs_i} +\gr(h_{\gs_i})-1
&\mbox{ for } i\in [t-1].\\
	  h_{\go_i} +\gr(h_{\go_i})&
 \mbox{ for } \;\; t\le i \le t+s-2.
\end{array} \right. \ee
Then let $H_J=\prod_{i \in r(J)}H_{i}$, and note that by \eqref{6t} and \eqref{9t}, we have
by
\be \label{3s} H_{i+r-1}(\gl) = 
 A_{i+r-1} \mbox{ for } i\in [t+s-2].
\ee
By Theorem \ref{7t}, ${\stackrel{\longrightarrow }{{\rm det}}}\;\A^{r,s}(\gl)$
is the evaluation of the \v Sapovalov element $\gth_\gc$ at $\gl \in \cH_\gc$.
We use the same procedure as before to  obtain the complete expansion of
${\stackrel{\longrightarrow}{{\rm det}}}\;\A^{r,s}(\gl)$ and thus the
\v Sapovalov element $\gth_\gc \in U(\fgl(t,s)) \subseteq U(\fgl(m,n)).$
\bt\label{bsb2} The \v Sapovalov element $\gth_\gc$ for the root $\gc=\gep_r - \gd_s$ in $U(\fgl(m,n))$  is given by
\be\label{2s1}\gth_\gc= \sum_{J\subseteq \I} f_J H_J.\ee
\et
\bpf   This follows from the proof of Theorem \ref{bsb1}. \epf
\bc \label{m2}For $\gl\in \cH_\gc$, $\gth_\gc(\gl)={\stackrel{\longrightarrow }{{\rm det}}}\;\A^{r,s}(\gl)$. \ec
\bpf By Equation \eqref{3s},
\[\sum_{J\subseteq \I} f_J H_J(\gl)  = \sum_{J\subseteq \I} f_J \prod_{i \in r(J)}A_{i}
\] and the argument in Theorem \ref{hea} shows that in the complete expansion of ${\stackrel{\longrightarrow }{{\rm det}}}\;\A^{r,s}(\gl)$, the coefficient of $f_J$ is $\prod_{i \in r(J)}A_{i}$.
\epf
\subsection{More Expansions for \v Sapovalov elements.}\label{exps}
We consider two new orders on root vectors, in one the unique odd root vector appears last in each product $f_J$.  In the other they appear first.  Applications of each are given in Subsection \ref{mrs}.
\subsubsection{Putting odd root vectors last.}\label{1exps}

With $\I$ as in \eqref{2w}, we give another expansion of the \v Sapovalov element $\gth_\gc $  using a different order on positive order roots. With this new order for each ${I\subseteq \I}$, the unique odd root vector appears last. Because of this  requirement, the definition of  $f_I$ is not as easily described using an ordering on $I$.    Thus consider $I$ as an unordered set 
\be \label{ords} I= \{m+s, i_1,i_2 \ldots, i_{g},    j_{1}, \ldots ,j_{h}\}\ee
 where $m+s > i_1 >i_2 \ldots > i_{g}\ge m+1$ and $r=  j_{1}< \ldots <j_{h}\le m$. Then  we define
\be \label{987t}f_I= e_{m+s, i_1}e_{i_1,i_2} \ldots,e_{ i_{g-1}, i_{g}} e_{j_2, j_1}\ldots,e_{ j_{h}, j_{h-1}}e_{i_{g}, j_h}.\ee 
 Note that the odd root vector $e_{i_{g}, j_h} $ is the last entry in $ f_I$.   
The definition of $H_I$ is also changed.  
Instead of \eqref{9t}, we define $H_{i+r-1}\in U(\fh)$ by
\be \label{98t} H_{i+r-1} = \left\{ \begin{array}
  {cc} h_{\gs_i} +\gr(h_{\gs_i})
&\mbox{ for } i\in [t-1]\\
	h_{\go_i} +\gr(h_{\go_i})&
 \mbox{ for } \;\; t\le i \le t+s-2.
\end{array} \right. \ee 
Recall the definition of  $A_i $ from \eqref{6t} and set $C_i=A_i +1$ for all $i$. Then
\be \label{98ts} H_{i+r-1}(\gl) = \left\{ \begin{array}
  {cc} C_{i+r-1}
&
\mbox{ for } i\in [t-1]\\
A_{i+r-1}
&
 \mbox{ for } \;\; t\le i \le t+s-2.
\end{array} \right. \ee 
Define $r(I)$ as usual, then set $H_I=\prod_{i \in r(I)}H_{i}$. Then we have 
\bt\label{bsb3} With the above notation, the \v Sapovalov element $\gth_\gc$ for the root $\gc=\gep_r - \gd_s$ in $U(\fgl(m,n))$  is given by
\be\label{2s}\gth_\gc= \sum_{I\subseteq \I} f_I H_I.\ee
\et
\bpf 
Recall the matrix  $\A^{r,s}(\gl) $ defined in 
\eqref{gbv}. The entries in  $\A^{r,s}(\gl) $ belong to the subalgebra $\fg'\cong \fgl(t,s)$. Note that the entries in the last $t$ columns and first $s$ rows of  $\A^{r,s}(\gl) $  are precisely the odd negative root vectors of $\fg'.$
We use the following {\it cofactor procedure}.
Begin with a  cofactor expansion of $\A^{r,s}(\gl) $   down the first column.   
Then use a cofactor expansion of the resulting cofactors down their first columns. 
Repeat this  a total of $s-1$ times. Note that every time we use cofactor expansion on a minor, the two resulting cofactors have  one fewer row containing odd elements.  Thus at the end of the procedure, all the remaining minors contain only odd root vectors in the first row, and all root vectors in this row are odd.  Next we use Proposition  \ref{2aa} to replace each minor by a matrix containing only odd vectors in the last column. Then the result follows  easily after some bookkeeping. Let $\K = \{K\subseteq \{m+1,\ldots ,m+s\}| m+s\in K\}$. For $K\in \K$, put the entries  of $K$ is descending order and define
$f_K$  as in \eqref{ddd}. 
Denote the complement of $K$ in $\{m+1,\ldots ,m+s\}$ 
by $\bar K$, and set $s(K) =\{p-1|p \in \bar K \}$. 
The smallest element of $K$ is $m+j$ for some $j\in[s] $, and we set  $j(K)=j$ in this situation.
Next we describe the minors that can arise at the end of the procedure.  First 
let 
\be \label{gwv} \F(\gl) = \left[ {\begin{array}{cccccc}  
-A_{m-1} &e_{m,m-1}&e_{m,m-2}& \hdots& e_{m,r+1}& e_{m,r}\\
0&-A_{m-2} &e_{m-1,m-2}&\hdots&  \hdots& e_{m-1,r}\\
\vdots& \vdots  &\ddots & \ddots &\vdots& \vdots\\
\vdots& \vdots  & \hdots& -A_{r+1}& e_{r+2,r+1} &   e_{r+2,r} \\
0&  0  & 0&0&  -A_r& e_{r+1,r} \\
\end{array}}
 \right],
\ee  and note that $\F(\gl)$ fits into the bottom right corner of $\A^{r,s}(\gl) $.
Now let $\F^{(j)}(\gl)$ be the matrix obained from $\F(\gl)$ obtained by adjoining, as the first row the vector
\be \label{gaav} (e_{m+j,m},e_{m+j,m-1},\hdots, \hdots, e_{m+j,r+1}, e_{m+j,r}).\ee
Note that all entries in \eqref{gaav} are odd, while all other entries in $\F^{(j)}(\gl)$ not on the 
subdiagonal belong to $\fgl(m)\op 0.$ 
Then the result of the   cofactor procedure is the following.
\bl If $\gl \in  \cH_\gc$, then 
\[\gth_{\gc}(\gl)= \sum_{K \in \K} f_K \prod_{k \in s(K)}  A_k \;
{\stackrel{\longrightarrow }{{\rm det}}}(\F^{(j(K))}(\gl))  \]

\el \bpf In each step of the cofactor procedure  we choose either an element from $\fn^-$ or a subdiagonal element of $\A^{r,s}(\gl) $. In column $i$, the entries in  $\fn^-$ all have the  form  $e_{*,m+s-i }$ and the subdiagonal element is $-A_{m+s-i -1}$. We choose an entry from  $\fn^-$ from this column iff ${m+s-i}\in K$.  Otherwise we choose $-A_{m+s-i -1}$, which agrees with ${m+s-i -1\in s(K)}$. Now in the complete expansion of 
$\A^{r,s}(\gl) $, 
$f_K$ is the initial part of various products  $f_I$ with $I\in \I$. If $j=j(K)$, then the last factor in $f_K$ has the form
$e_{*,m+j }$. Then to get a valid product the next term should be  $e_{m+j,* }$ which means that the minor associated to $K$ at the end of  the cofactor procedure is
 $\F^{(j(K))}(\gl)$. 
 \epf \noi
Now  let $ \G^{(j)}(\gl)$ be the matrix  
\be \label{gvv} \G^{(j)}(\gl) = \left[ {\begin{array}{cccccc}  
e_{r+1,r}&e_{r+2,r}&e_{r+3,r}& \hdots & e_{m,r}& e_{m+j,r}\\
-C_{r} &e_{r+2,r+1}&e_{r+2,r+1}& \hdots& e_{m,r+1}& e_{m+j,r+1}\\
0&-C_{r+1} &e_{r+3,r+2}&\hdots&  \hdots& e_{m+j,r+2}\\
\vdots& \vdots  &\ddots & \ddots &\vdots& \vdots\\
\vdots& \vdots  & \hdots& -C_{m-2}& e_{m,m-1} &   e_{m+j,m-1} \\
0&  0  & 0&0&  -C_{m-1}& e_{m+j,m} \\
\end{array}}
 \right].
\ee Then by Proposition  \ref{2aa},
$
{\stackrel{\longrightarrow }{{\rm det}}}(\F^{(j)}(\gl))  = {\stackrel{\longrightarrow }{{\rm det}}}(\G^{(j)}(\gl)) $. Hence
\[\gth_{\gc}(\gl)= \sum_{K \in \K} f_K \prod_{k \in s(K)}  A_k \;
{\stackrel{\longrightarrow }{{\rm det}}}(\G^{(j(K))}(\gl)).  \]
 Let $\J(j)=\{J \subseteq \{r, r+1, \ldots, m, m+j \} | r, m+j \in J\}.$  It is clear that in the complete expansion of the determinant ${\stackrel{\longrightarrow }{{\rm det}}}(\G^{(j)}(\gl)) $, we 
encounter only terms $f_J$, for some $J\in \J(j)$ ordered as in the last part of \eqref{987t}, and that all such $J$ occur. To complete the proof we use a slightly different cofactor procedure than before. We start with a cofactor expansion along the last row, which we call row 1, and continue upwards, numbering rows from bottom to top. In row $i$ the entries from $\fn^-$ have the form $e_{*, m+1-i}$ and the subdiagonal entry is $-C_{m-i}$.  Now $m+1-i \in J$ iff we choose an element from $\fn^-$ from this row.  Otherwise we choose 
$-C_{m-i}$, which agrees with ${m-i}\in r(J)$. Now define $t(J)$  using the complement of $J$ in $\J(j)$. Then we have 
\[  
{\stackrel{\longrightarrow }{{\rm det}}}(\G^{(j)}(\gl)) 
= \sum_{J \in \J(j)} f_J \prod_{i \in t(J)} 
C _i. \] 
Thus 

\[\gth_{\gc}(\gl)= \sum_{K \in \K} f_K \sum_{J \in \J(j(K))} f_J 
\prod_{i \in s(K)}  A_i 
\prod_{i \in t(J)} 
C _i. \] 
It is clear that for any $K\in \K$ and $J\in \J(j)$ with $j=j(K)$, we have $f_Kf_J \in U(\fn^-)^{-\gc}$ and the terms are correctly ordered so $f_Kf_J  =f_I$ where
$I=K\cup J \in \I$. We have 
$r(I)=s(K) \cup  t(J)$. 
Now to complete the proof of Theorem \ref{bsb3}, we need to show that 
$$ \prod_{i \in r(I)}H_{i}(\gl)=\prod_{i \in s(K)}  A_i 
\prod_{i \in t(J)}C_i.$$  This follows from \eqref{98ts}.
\epf
\subsubsection{Putting odd root vectors first.}\label{2exps}
  If $I$ is the set in \eqref{ords} define 

\be \label{987st}f_I=e_{i_{g}, j_h}  e_{j_2, j_1}\ldots,e_{ j_{h}, j_{h-1}}e_{m+s, i_1}
e_{i_1,i_2} \ldots,e_{ i_{g-1}, i_{g}}.\ee
 Next 
instead of \eqref{98t}, we define $H_{i+r-1}\in U(\fh)$ by
\be \label{98st} H_{i+r-1} = \left\{ \begin{array}
  {cc} h_{\gs_i} +\gr(h_{\gs_i})-1
&\mbox{ for } i\in [t-1]\\
	h_{\go_i} +\gr(h_{\go_i})+1&
 \mbox{ for } \;\; t\le i \le t+s-2.
\end{array} \right. \ee 
Note that in place of \eqref{98ts}, we now have  
\be \label{58ts} H_{i+r-1}(\gl) = \left\{ \begin{array}
  {cc} A_{i+r-1}
&
\mbox{ for } i\in [t-1]\\
C_{i+r-1}
&
 \mbox{ for } \;\; t\le i \le t+s-2.
\end{array} \right. \ee 

\bt\label{bab3} With the above notation, the \v Sapovalov element $\gth_\gc$ for the root $\gc=\gep_r - \gd_s$ in $U(\fgl(m,n))$  is given by
\be\label{22s}\gth_\gc= \sum_{I\subseteq \I} f_I H_I.\ee
\et
\bpf This is similar to the proof of Theorem \ref{bsb3}, so we give fewer details. We use a cofactor expansion
of $ \A^{r,s}(\gl)$ along the last $t-1$ rows working from the bottom up. 
The subdiagonal entries that emerge in these steps are the $A_{i+r-1}$ 
 for $ i\in [t-1]$.
  At the end of these steps, the unexpanded cofactors have last column consisting of odd root vectors and no other odd root vectors.   Then  we use Proposition  \ref{2aa} to replace each of these cofactors with a matrix containing only odd vectors in the first row, and subdiagonal entries 
$C_{i+r-1}$ for $t \le i \le t+s-2.$ The proof concludes as before.
\epf

\subsection{Partial Expansions of \v Sapovalov elements.}\label{expa}
We obtain two expansions of the \v Sapovalov element for the root $\gc=\gep_r - \gd_s$ in \eqref{2s} under some additional assumptions. These expansions correspond to partial expansions of the determinant
${\stackrel{\longrightarrow }{{\rm det}}}\;\A^{r,s}(\gl)$, obtained by separating out the terms containing certain subdiagonal entries, whence the title of this subsection. The main results are Theorems \ref{stc} and \ref{bsb}.  Motivation for these results is given in Subsection \ref{mrs}.
The additional assumptions are\\ \\
{\bf Case 1} We assume that $(\gl+\gr,\gc) = (\gl+\gr, \gc') =0$, where $\gc' = s_\alpha\gc = \gc -\ga$ for some positive root $\ga$. 
This means that $(\gl + \gr,\ga^\vee ) =0$ and by  \cite{M2} Theorem 5.9 
 we have, in this situation $\theta_{\gc}v_{\gl}= \theta_{\ga,1} \theta_{\gc'} v_{\gl}$.
Theorem \ref{stc} below (the main result in Case 1) could be considered as a refinement of this equation when $\gl$ is replaced by
${{\widetilde{\lambda}}}=\gl+\ttT\xi$. 
 \\ \\
{\bf Case 2}
We suppose $\ga_1, \ga_2$ (resp. $\gc_1, \gc_2$) are distinct positive even (resp. odd) roots	
such that
\be \label{bla} \ga_1 +\gc_1 +\ga_2 \;=\gc_2.\ee
This means that 
$\ga_1^\vee \equiv\ga_2^\vee \mod \Z\gc_1 +\Z\gc_2$.  
If $ \gc_1, \gc_2 \in X$  there are significant consequences for the structure of the modules $M^{X}({{\lambda} })$ from \cite{M2} Theorem 1.10, see \ref{hut}, and factors of the corresponding  \v Sapovalov determinant $\det F^X_{\eta}$. 
\\ \\
We disregard  the notation for $\ga_i, \gc_i$ introduced at the start of Section \ref{1s.8}. Then for Case 1 we assume that 
\be \label{ha} \gc =\gep_r- \gd_{s}, \ga = \gep_r-  \gep_{\ell }  \mbox{ and } \gc'=s_\ga\gc=\gep_\ell- \gd_{s}.\ee In Case 2, we may assume that 
\be \label{hz} \ga_1= \gep_r-  \gep_{\ell }, \ga_2 = \gd_k- \gd_{s} \mbox{ and }  \gc_1 =\gep_{\ell }- \gd_k, \gc_2 =\gep_r- \gd_{s}.\ee

\subsubsection{Case 1.} \label{case1}
We assume that Case 1 and \eqref{ha} hold. 
Before we state the main result there is another issue to deal with. We need to use a different ordering on positive roots.
To do this we introduce the matrix

\be \label{grv} \B^{r,s}(\gl) = \left[ {\begin{array}{cccccc}
e_{r+1,r}&e_{r+2,r} &\hdots &\hdots &e_{m+s-1,r}  & e_{m+s,r} \\
 -C_{r} &e_{r+2,r+1} &\hdots &\hdots&e_{m+s-1,r+1} &  e_{m+s,r+1} \\
 \vdots  & \vdots & \ddots&\hdots &\hdots& \hdots \\
 \vdots  & \vdots & \ddots&\hdots &\hdots& \hdots \\
0&0& \hdots &  -C_{m+s-1} & e_{m+s-1,m+s-2} &   e_{m+s,m+s-2} \\
0&0&\hdots&  \hdots  & -C_{m+s-2}& e_{m+s,m+s-1} \\
\end{array}}
 \right],
\ee 
By Proposition \ref{2aa} and \eqref{gbv}, we have
\be \label{sa2}{\stackrel{\longrightarrow }{{\rm det}}}\;\A^{r,s}(\gl) ={\stackrel{\longrightarrow }{{\rm det}}}\;\B^{r,s}(\gl).\ee
 At this point we bring $\ga = \gep_r - \gep_\ell = \gs_{\ell-r}$ and $\ttT= H_{\ell-1}$ into the story by noting that $(\ell-r) +(r-1)=\ell-1$, so by Equations \eqref{6t} and \eqref{98ts},
 we have
\be \label{35b} \ttT(\gl) = C_{\ell-1}=(\gl+\gr,\ga_{1}). \ee
With $\I$ as in \eqref{2w}, we define
\[ \I_{\ell } = \{I\subseteq \I| \ell -1 \in r(I)\},\;\;\I^{\ell } = \{I\subseteq \I| \ell-1\notin r(I)\}.\]
Note that  $I\subseteq \I^\ell $ if and only if $\ell \in I$. Also
$\ttT= H_{\ell -1}$ is a factor of $H_J$ iff 
 ${J  \subseteq \I_{\ell}}$. 
For ${J  \subseteq \I_{\ell}}$, define $H'_J =
H_J/\ttT$.

\bt \label{stc} If $\gth_{\gc}$ is as in Theorem \ref{bsb2}, then for $\gl\in\cH_{\gc},$ we have
\be\label{mice3} 
\gth_{\gc}v_\gl=
(\gth_{\ga,1}\gth_{\gc'}+\sum_{J  \subseteq \I_{\ell}}f_J H'_J \ttT)v_\gl
.\ee
\et
\bpf
We have \be\label{miice}
\gth_{\gc}=
\sum_{J  \subseteq \I^{\ell}}f_J H_J+\sum_{J  \subseteq \I_{\ell}}f_J H'_J \ttT
.\ee
\noi
It is useful to think about this proof in terms of non-commutative determinants. Thus we set
$B= \B^{r,s}(\gl) $ as in \eqref{grv},
and  recall that $\ttT(\gl) = C_{\ell-1}$ by \eqref{35b}.
By Corollary \ref{m2} and \eqref{sa2} we have 
$\gth_{\gc}(\gl)={\stackrel{\longrightarrow}{{\rm det}}} B$. 
From \eqref{sam2},
\be \label{sim}{\stackrel{\longrightarrow}{{\rm det}}} B =
{\stackrel{\longrightarrow}{{\rm det}}}B' + \ttT(\gl)  {\stackrel{\longrightarrow}{{\rm det}}}B'',\ee
where  $B'$ and  $B''$ are obtained from  $B$ by setting $\ttT(\gl) =0,$
and by deleting the row and column containing $\ttT(\gl)$ respectively. Hence $B''$ is the cofactor of entry  $\ttT(\gl)$ in $B$.
Also any term in  the expression \eqref{lr} for ${\stackrel{\longrightarrow }{{\rm det}}}(B)$ which contains a factor  of the form $e_{*, \ell}$ or $e_{\ell,*},$ 
 cannot occur in ${\stackrel{\longrightarrow}{{\rm det}}} B''$.
This means that  ${\stackrel{\longrightarrow}{{\rm det}}}B''$ belong to $U(\overline{\fg})$ where
$\overline{\fg}$ the  subalgebra isomorphic to $\fgl(m-1,n)$ and
with rows and columns indexed by the set
$\I(\ell)=\{r, \ldots, \widehat{\ell},\ldots, m+s\}$.
As mentioned above $J\subseteq \I_\ell $ if and only if $\ell \notin J$. 
Thus the  second term on the right of  \eqref{sim} equals 
the evaluation of the second term on the right of \eqref{miice} at  $\gl$. 
Note that when $\ttT(\gl)=0$, $\gth_\gc v_{{\lambda}} =
\gth_{\ga,1}\gth_{\gc'}v_{{\lambda}},$ as already mentioned, so the proof is complete.
\epf \noi
\subsubsection{Case 2.} \label{case2} 
We assume that Case 1 and \eqref{hz} hold. 
The proof of our main result depends on Theorem \ref{bsb3}. Hence we  define $f_I, A_i, H_i$ as in  \eqref{987t}, \eqref{6t} and \eqref{98t} respectively. The relation between these elements is given by \eqref{98ts}. 
Recall $\ga_1 = \gep_r - \gep_\ell = \gs_{\ell-r}$ and $\ga_2= \gd_k -\gd_s $ 
and we have 
\be\label{1ff} A_{\ell-1} = (\gl +\gr,\ga_{1})-1 \mbox{ and } A_{m+k - 1} = (\gl+\gr,\ga_{2}).\ee
It may help to observe that in \eqref{67t} $\go_i= \gd_k -\gd_s  =\ga_2$, when $i=k+t-1$, and thus $i+r-1 = {m+k - 1}$ in \eqref{98ts}.
\\ \\
We treat 
\be \label{w12} \ttT= H_{\ell-1} \mbox { and }\ttS=-H_{m+k-1}\ee
 as indeterminates which can be evaluated on any $\gl\in \cH_\gc$. By Equation \eqref{98ts}, we have
\be \label{3b} \ttT(\gl) = C_{\ell-1}=(\gl+\gr,\ga^\vee_{1}) \mbox{ and } \ttS(\gl) =- A_{m+k - 1} = (\gl+\gr,\ga^\vee_{2}). \ee
For each subset $J$ in \eqref{2s} we are interested in when $\ttS$ or $\ttT,$ from \eqref{w12},  or both are factors of $H_J$.
Note that
\by \label{help} \ttT  \mbox{ is a factor of } H_J \mbox{ iff }
\ell -1 \in r(J), \mbox{ iff } J \in \I(\ga_2) \mbox{ or }\I(\ga_1, \ga_2)\nn\\
\ttS  \mbox{ is a factor of } H_J \mbox{ iff }
m+k-1 \in r(J) \mbox{ iff }
J \in \I(\ga_1) \mbox{ or }\I(\ga_1, \ga_2).
 \ey
where
\by \label{79t} \I({\emptyset}) = \{J\subseteq \I| \ell-1, m+k-1 \notin r(J)\}, \;
\I(\ga_1) = \{J\subseteq \I| \ell-1 \notin r(J), m+k-1 \in r(J)\}, \;\;\;\;\\
\I(\ga_2) = \{J\subseteq \I| m+k-1 \notin r(J), \ell-1 \in r(J)\},\;\;\;\nn
\I(\ga_1, \ga_2) = \{J\subseteq \I| \ell-1, m+k-1 \in r(J)\}.\;\;\;\ey
\noi The notation is inspired by Equation  \eqref{3b}. 
Alternatively
\by \label{71} \I({\emptyset}) = \{J\subseteq \I| \ell, m+k \in J\}, \;
\I(\ga_1) = \{J\subseteq \I| \ell \in J, m+k \notin J\}, \;\\
\I(\ga_2) = \{J\subseteq \I| m+k \in J, \ell \notin J\}.\;\;\;
\I(\ga_1, \ga_2) = \{J\subseteq \I| \ell, m+k \notin J\}.\;\;\;\nn\ey
Observe that $\I$ is the disjoint union of the above 4 sets.
Next for $J$ a subset of $\I(\ga_1), \I(\ga_2), \I(\ga_1, \ga_2)$ elements
$H_J^{(1)},\;  H_J^{(2)},\; H_J^{(3)}$ respectively by
\be \label{711} H_J^{(1)} =  H_J/\ttS ,\; H_J^{(2)}=  H_J/\ttT  ,\; H_J^{(3)} =  H_J/\ttS \ttT.
\ee By \eqref{help} these are all elements of  $U(\fh)$.
\bt\label{bsb} The \v Sapovalov element $\gth_\gc$ for the root $\gc=\gep_r - \gd_s$ in satisfies
\be\label{mice}
\gth_{\gc}=\gth_{\ga_1}\gth_{\ga_2}\gth_{\gc_1} -
\sum_{J \in\I{(\ga_1)}} f_J H_J^{(1)}\ttS 
+\sum_{J  \in  \I{(\ga_2)}}f_J H_J^{(2)}\ttT
 \;- \sum_{J \in   \I(\ga_1, \ga_2)}f_J H_J^{(3)}\ttS\ttT
.\ee
\et
\bpf 
We break the sum in Theorem \ref{2s} into 4 pieces corresponding to the disjoint union of the sets in \eqref{71}. Then  we have to show that $\gth_{\ga_1}\gth_{\ga_2}\gth_{\gc_1} = \sum_{J  \in \I({\emptyset})}f_J H_J$.
 Recall that $\ttT(\gl)= C_{\ell-1}$ and $\ttS(\gl)=-A_{m+k-1}.$ We first consider what happens when these are zero. Since $C_{\ell-1}=0$, $\B^{r,s}(\gl)$ is block upper triangular, and this gives
${\stackrel{\longrightarrow }{{\rm det}}}\;\B^{r,s}(\gl) = \gth_{\ga_1}{\stackrel{\longrightarrow }{{\rm det}}}\;\B^{\ell,s}(\gl) =  \gth_{\ga_1}{\stackrel{\longrightarrow }{{\rm det}}}\;\A^{\ell,s}(\gl).$ Then $A_{m+k-1}=0$ implies that $\A^{\ell,s}(\gl)$ is block upper triangular, and we deduce that $\gth_{\gc}=\gth_{\ga_1}\gth_{\ga_2}\gth_{\gc_1}$. Now $J \in \I({\emptyset})$ means that
$\ell-1, m+k-1 \notin r(J)$, and thus neither of $\ttS$ or $\ttT$ is a factor of $H_J$.  So the sum $\sum_{J  \in \I({\emptyset})}f_J H_J$ is independent of  $\ttT(\gl)$ and $\ttS(\gl)$.
It follows that
$\gth_{\ga_1}\gth_{\ga_2}\gth_{\gc_1} = \sum_{J  \in \I({\emptyset})}f_J H_J$ as required.
\epf
\br \label{231} { \rm From the proof of Theorem \ref{stc}, that the  term
$\sum_{J  \in \I_{\ell}}f_J H'_J$ in is the  \v Sapovalov element for the longest root $\bar \gc$ of the subalgebra $\overline{\fg}$ of  $\fg$.
To see this observe that with ${\stackrel{\longrightarrow}{{\rm det}}}B''$ as in the proof, we have $\gth_{\bar \gc}(\gl) ={\stackrel{\longrightarrow}{{\rm det}}}B''$.
 We can make a similar remark about the last three sums on the right of \eqref{mice}.  To do this we need the index sets,  cf. \eqref{71},
\by \label{2wa}
\I_{[\ga_1]} &=&  \{r, r+1,\ldots,\widehat{\ell},\ldots ,m+s\},\nn\\
\I_{[\ga_2]} &=& \{r, r+1,\ldots,\widehat{m+k},\ldots m+s\}\\
\I_{[\ga_1, \ga_2]}
&=& \{r, r+1,\ldots,\widehat{\ell},\ldots \widehat{m+k}, m+s\},\nn
.\ey
The set $\I_{[\ga_1]}$ indexes the rows and columns of the
subalgebra
$\overline{\fg} \cong \fgl(m-1,n)$ as in Case 1. Define  subalgebras
$\underline{\fg}\cong\fgl(m,n-1)$
and
$\underline{\overline{\fg}} \cong\fgl(m-1,n-1)$, similarly, using the sets
$\I_{[\ga_2]} $ and $\I_{[\ga_1, \ga_2]}$ respectively 
to index the rows and columns. Then in Equation \eqref{mice} from Theorem \ref{bsb}, the coefficients of $\ttS, \ttT$ and $\ttS \ttT$ are the \v Sapovalov elements for the longest odd root of the subalgebras $\underline{\fg}, \overline{\fg}$  and
$\underline{\overline{\fg}}$
respectively. More details on this approach were given in \cite{M2}. Since we do not need this extra information about Theorems \ref{stc} or \ref{bsb}, we give no further details here.
}\er

\bexa \label{22ex}{\rm Let $\fg= \fgl(2,2)$ and using the notation of \eqref{heb},  set $ \ga = \ga_1$ and $ \gc = \gc_1$.  We also need 
\[
e_{-\ga -\gb} = e_{31}, \quad  
e_{-\gb-\gc} = e_{42}, \quad e_{-\ga-\gb-\gc} = e_{41}. 
\]
 We find the \v Sapovalov elements for the roots $\ga + \gb,
{\gb+\gc}$ and ${\ga+ \gb+\gc}$ using Theorem \ref{bsb1}.
\by \label{a+b}\gth_{\ga +\gb}&=&  e_{-\gb}e_{-\ga}+e_{-\ga -\gb}h_\ga \\
&=&  e_{-\ga}e_{-\gb}+e_{-\ga -\gb}(h_\ga +1). \nn\ey
\by \label{c+b} \gth_{\gb+\gc}&=& e_{-\gc}e_{-\gb}+e_{-\gb-\gc}(h_\gc-1)\\
&=& e_{-\gb}e_{-\gc}+e_{-\gb-\gc}h_\gc\nn
. \ey
We give four  expressions for $\gth_{\ga+\gb+\gc}$. The first using Theorem \ref{bsb1}, the second using Theorem \ref{bsb3} or  Theorem \ref{bsb}, the third using 
 Theorem \ref{bab3} and the fourth by expanding the determinant of the matrix in \eqref{grv}.
\by \label{97st} \gth_{\ga+\gb+\gc}
   & =&  e_{-\gc}e_{-\gb}e_{-\ga}+
e_{-\gc}e_{-\ga -\gb}h_\ga 
+e_{-\gb-\gc}e_{-\ga}(h_\gc-1) +e_{-\ga-\gb-\gc}h_\ga (h_\gc-1)\nn
\\& =&	e_{-\gc}e_{-\ga}e_{-\gb}+
e_{-\gc}e_{-\ga -\gb}(h_\ga +1)
+e_{-\ga}e_{-\gb-\gc}(h_\gc-1) +e_{-\ga-\gb-\gc}(h_\ga +1)(h_\gc-1)\nn\\
& =&	 e_{-\gb}e_{-\gc}e_{-\ga}+
e_{-\ga -\gb}e_{-\gc}h_\ga 
+e_{-\gb-\gc}e_{-\ga}h_\gc +e_{-\ga-\gb-\gc}h_\ga h_\gc\\
& =& e_{-\ga}e_{-\gb}e_{-\gc}+
e_{-\ga -\gb}e_{-\gc}(h_\ga +1)
+e_{-\ga}e_{-\gb-\gc}h_\gc +e_{-\ga-\gb-\gc}(h_\ga +1)h_\gc\nn
. \ey
Using the relations in $U(\fg)$, the right sides of  the above expressions are easily seen to be equal.

}\eexa

\subsection{Motivation for the Expansions.} \label{mrs} 
Let $X$ be an orthogonal set of positive isotropic roots. Suppose 
$(\gl + \gr, \gc) =0$ for all $\gc \in X.$ In \cite{M2}, see also \cite{M4}
we constructed some highest weight modules $M^X(\lambda)$ with highest weight $\gl$, and character $ \tte^{\gl} p_X$, where $p_X$ is a partition function that counts partitions not involving roots in $X$. The definition of these modules is as follows.
Let $T$ be an indeterminate and set  $A=\ttk[T], B=\ttk(T)$.  For $R=A$ or $B$, set $U(\fg)_R= U(\fg)\ot_\ttk R.$  
Now choose $\xi\in \fh^*$ subject to certain conditions which are spelled out in the two cases below, but otherwise generic, 
 set ${{\widetilde{\lambda}}}=\gl+T \xi$, and consider the Verma module $M({\widetilde{\lambda}})_B$ over  $U(\fg)_B$ with highest weight ${\widetilde{\lambda}}$. Then 
set 
\be \label{tuv} M^X({\widetilde{\lambda}})_B = M({\widetilde{\lambda}})_B/ \sum_{\gc \in X} U(\fg)_B \gth_\gc v_{{\widetilde{\lambda}}}.\ee
Next let $M^{X}({\widetilde{\lambda} })_{A}$ be the
$U(\fg)_A$-submodule of
$M^{X}({\widetilde{\lambda} })_{B}$ generated by the highest weight vector and define 
\be  \label{hut} M^{X}(\lambda) = M^{X}({\widetilde{\lambda} })_{A}/TM^{X}({\widetilde{\lambda} })_{A}.\ee
In \cite{M2} Section 11, we evaluate the \v Sapovalov determinant $\det F^X_{\eta}$ and give a Jantzen sum formula  for these modules. A difficulty which does not arise in the case of Verma modules \cite{M} Theorem 10.2.5, is that there is no natural $A$-basis for the weight space  $M^{X}({\widetilde{\lambda} })_{A}^{\widetilde{\lambda} - \eta}$ 
indexed by partitions.  To surmount this difficulty we use related  determinant  $\det G^X_{\eta}$ of a matix $ G^X_{\eta}$ with rows and columns indexed by elements $e_{-\pi} v_{\widetilde{\lambda} } $, with $\pi$ a partition of  $\eta$  as in the classical  case.  The leading term of $\det G^X_{\eta}$ is easy to compute. In 
 Theorem  11.1 of \cite{M2} we determine the relationship between the leading terms of  the two determinants  $\det G^X_{\eta}$ and $\det F^X_{\eta}$, and evaluate the latter.
The comparison of leading terms relies on Equations \eqref{msd} and \eqref{maace} below, which are used repeatedly to improve the $A$-basis  used to calculate $\det G^X_{\eta}$. 
Below we refer to Cases 1 and 2 from the previous Subsection. 
\\
\\Case 1:  Let $T =\ttT$ be as in \eqref{35b}. 
Suppose  $(\gl+\gr,\gc) = (\gl+\gr, \gc') =0$, where $\gc' = s_\alpha\gc = \gc -\ga$ for some positive root $\ga$. 
Then  $(\gl + \gr,\ga^\vee ) =0$.
Choose $\xi\in \fh^*$ so that
$(\xi,\gc)= (\xi,\gc')=0=(\xi,\ga^\vee) $ and with $\xi$ otherwise generic. 
Then set ${{\widetilde{\lambda}}}=\gl+T \xi$.
Using the notation of \eqref{mice3}, consider the  elements of $M({\widetilde{\lambda}})$  given by
$$p= \gth_{\ga,1}\gth_{\gc'}v_{\widetilde{\lambda}}\mbox{ and } q=
\sum_{J  \subseteq \I_{\ell}}f_J H'_J   v_{\widetilde{\lambda}},$$ 
we have  the relation $\gth_\gc v_{\widetilde{\lambda}}= p+qT$ in $M({\widetilde{\lambda}})_A$. 
Thus if $\gc\in X$, we have  in  $M^{X}({\widetilde{\lambda} })_{A}$.
\be \label{msd} p+qT =0.\ee
This is a key step in the proof of Theorem 11.4 of \cite{M2}.
\\ \\
Case 2: Recall by \eqref{3b} that  $\ttT(\gl) =(\gl+\gr,\ga^\vee_{1}) \mbox{ and } \ttS(\gl)  = (\gl+\gr,\ga^\vee_{2})$. Set $\ttS=\ttT = T+1$.  Thus if $T=0$, we have $(\gl+\gr, \alpha^\vee_1) = (\gl+\gr, \alpha^\vee_2)=1$ which is the assumption of Theorem 11.8 of \cite{M2}.
Assume  that $\gl\in \cH_{\gc_1} \cap \cH_{\gc_2}$.
Now choose $\xi\in \fh^*$ so that
$(\xi,\gc_1)= (\xi,\gc_2)=0$, and $(\xi,\ga^\vee_1) = (\xi,\ga^\vee_2) = 1$ and with $\xi$ otherwise generic. 
Then set ${{\widetilde{\lambda}}}=\gl+T \xi$, and consider the Verma module $M({\widetilde{\lambda}})$ over $B$ with highest weight ${\widetilde{\lambda}}$. 
Let  $M$ be the factor module of $M({\widetilde{\lambda}})_A$ defined by setting $\gth_{\gc_1}v_{\widetilde{\lambda}} = \gth_{\gc_2}v_{\widetilde{\lambda}} =0.$
Then $(\widetilde{\lambda}+\gr, \ga_1^\vee) = (\widetilde{\lambda}+\gr, \ga_2^\vee) = T+1$.
Hence  in the notation of Theorem \ref{bsb} we have $\ttS=\ttT = T+1$.
Thus since $T+1$ is not a zero divisor in $B$, if  $\gth_{\gc_1}v_{\widetilde{\lambda}}=\gth_{\gc_2}v_{\widetilde{\lambda}} =0$, then \eqref{mice} yields the following relation in $M$ 
\be \label{maace}
0 = [\sum_{J \in  \I{(\ga_1)}} f_J H_J^{(1)}
-\sum_{J  \in  \I{(\ga_2)}} f_J H_J^{(2)}+
\sum_{J \in  \I{(\ga_1, \ga_2)}}f_J H_J^{(3)}(T+1)]v_{\widetilde{\lambda}}
\ee
This gives a generalization of Equation (B.20) from \cite{M2} and is a key step in the proof of Theorem 11.8 of \cite{M2}.
\\ \\
Finally Theorem \ref{bab3} can be used to give a new proof of a result of Serganova, \cite{S2} Theorem 5.5, see Theorem \ref{1shapel}. For the proof it is important that root vectors are ordered 
so that the unique odd root vector in each $f_J$ is the first factor. For the proof of Theorem \ref{1boy}  the unique odd root vector needs to be the last factor in each $f_J$.

\section{Powers of \v Sapovalov elements.} \label{pse}
\noi
Powers of \v Sapovalov elements are dealt with in detail in Section 5.1 of \cite{M2}.
The method, which works for all non-exceptional algebras, depends on Lemma \ref{1768}, which is the basis of the original constuction of \v Sapovalov elements from \cite{Sh}.  We briefly review the details.
First write the root $\eta$  in the form $\eta =w \gb$ with $\gb$ a simple root and $w\in W$.  Next  by induction on the length of  the $w$, we construct 
$\gth_\eta \in U(\fn^-)^{-\eta}$ such that $\gth_\eta v_\gl$ is a highest weight vector for $M(\gl)$ whenever $\gl\in \cH_\eta$, and the coefficient of $e_{-\pi^0}$ in $\gth_\eta $ equals 1, see \eqref{rbt}. Then we fix a (non-unique) lifting of $\gth_\eta$  to an element $\Gt_\eta\in U(\fb^-)^{-\eta}= U(\fn^-)^{-\eta} \ot U(\fh)$ such that $
\Gt_\eta v_\gl= \Gt_\eta(\gl) v_\gl =\gth_\eta v_\gl$,
for all $\gl \in {\mathcal H}_{\eta}$. 
\\ \\
 For an isotropic root $\eta$, if $\gl\in \cH_\eta$,  then also $\gl-\eta\in \cH_\eta$, and by induction on the length of  the $w$ we have
\bp \label{11.6} Let $\eta$ be an odd positive root of $\fgl(m,n),$ and define $\Gt_\eta$  as in Theorem \ref{bsb1}. If $\gl \in \cH_\gc,$ then $\Gt_\eta^2v_\gl = 0.$\ep

\bpf Combine Theorem 5.1 from \cite{M2} and Theorem \ref{bsb1}. \epf \noi
For a non-isotropic root $\eta$, the situation is more delicate, because it is no  longer true that $\gl\in \cH_{\eta, p}$ implies 
$\gl-\eta\in \cH_{\eta, p}$, and thus we need 
to evaluate $\Gt_{\eta, p}$   at a point that is not in 
$\cH_{\eta, p}$. 
To resolve this issue, a uniform inductive  constuction 
of $\Gt_\eta$  is given in \cite{M2} and  \cite{M4}. This depends on a specific choice of $\gb$ and a shortest length expression for $w$. 
 Then $\Gt_\eta^p$ is a \v Sapovalov elements for the pair $(\eta,p),$ by \cite{M2} Theorem 5.8. \\ \\
Needless to say, this method of proof is completely different from the direct approach used in most of this paper.  However the two approaches are brought together by 
the two proofs of Theorem \ref{hea}, the second of which uses the same choice of $\gb$ as in \cite{M2}.  From this we can deduce the following

\bp \label{11.5} Let $\eta$ be a positive root of $\fgl(m),$ and define $\Gt_\eta$  as in Theorem \ref{bb}. Then, for $p\ge 1$, 
$\Gt_\eta^p$ is the \v Sapovalov element for the pair $(\eta,p)$.  In other words if   $\gl\in \cH_{\eta,p}$, we have
\[\Gt_\eta^p(\gl) =\gth_{\eta,p}(\gl).\]
\ep.
\bpf This follows from the above discussion and Theorem \ref{bb}. \epf
\noi  The above result may be viewed as a version of \cite{CL} Theorem 2.7.

\section{{Survival of \v{S}}apovalov elements in factor modules.}\label{1surv}
 Let  $v_{\gl}$ be a highest weight vector
in a Verma module $M(\gl)$ with highest weight $\gl,$ and suppose $\gc$ is an odd root with
$(\gl+\gr, \gc) =0$. We are interested in the condition that the image of
$\theta_{\gc}v_{\gl}$ is non-zero in a factor module $M$ of $M(\gl)$. If this is the case, we might say that $\theta_{\gc}v_{\gl}$ {\it survives} in $M$.  Throughout this section we assume that $\fg = \fgl(m,n)$. \v{S}apovalov elements and Verma modules are defined using the distinguished Borel subalgebra.   
\subsection{Independence of \v{S}apovalov elements.}\label{ioe}
Given $\gl \in \fh^*$  set
\by \label{rue} 
B(\gl) &=& \{ \gc \in {\Delta}^+_{1} | (\lambda + \rho,
\gc) = 0 \} .\nn \ey
If $\gc,\gc'$ are  positive non-orthogonal isotropic roots, then
$\gc'=s_\ga \gc$ for some even root $\ga$.  Assume that  
\be \label{prn}    (\gc, \alpha^\vee) =1 \mbox{ and } (\gl+\gr, \alpha^\vee) =0\ee 
In this situation
the next result relates \v Sapovalov elements for $\gc, \gc'$ and $\ga$. 
\bt  \label{man}
Let  $v_{\gl}$ be a highest weight vector in $M(\gl)$ and set $\gc' = s_\alpha\gc.$  If $ (\gl+\gr, \alpha^\vee)
=0$ and $\gc'\in B(\gl)$ we have 
\be \label{pin}\theta_{\gc} v_{\gl}= \theta_{\ga,1} \theta_{\gc'} v_{\gl}.\ee\et
\bpf This is a special case of \cite{M2} Theorem 5.9 (a). Note that by \eqref{prn} and the hypothesis, $\gc\in B(\gl)$.\epf
\noi 
Define a ``Bruhat order" $\le$ on $B(\gl)$ as follows. First write  $\gc' \tl \gc$ if $\gc-\gc'$ is a  positive even root $\ga$, $ (\gl+\gr, \alpha^\vee)
=0$ and $s_\ga \gc=\gc'$.
 The relation $\le$ is the transitive extension of $\tl.$ Thus $\gc' \le \gc$ if there is a sequence of  elements $\gc=\gc_0, \gc_1, \ldots \gc_n = \gc'$  such that $\gc_i \tl \gc_{i-1}$ for $i\in [n]$.
 Then 
we say that $\gc$ is $\gl$-{\it minimal} if $\gc' \le \gc$  with $\gc' \in B(\gl)$ implies that $\gc' = \gc$.
For $\gc \in B(\gl)$ set $B(\gl)^{-\gc} = B(\gl)\backslash \{\gc\}$. We say $\gc$ is {\it independent
at } $\gl$ if $$\theta_{\gc}v_{\gl}\notin \sum_{\gc' \in B(\gl)^{-\gc} }U(\fg)\theta_{\gc'} v_{\gl}.$$
\bp \label{1pot}If $\gc' \tl \gc$ with $\gc' \in B(\gl)$, then $\theta_{\gc}v_{\gl}\in U(\fg)\theta_{\gc'} v_{\gl}.$ \ep
\bpf By  hypothesis \eqref{prn} holds. 
Thus the   result follows from Theorem \ref{man}.\epf
\noi By the  Proposition,  if we are interested in the independence of the \v{S}apovalov elements $\gth_\gc$ for distinct isotropic roots, it suffices to  study only $\gl$-minimal roots $\gc$.
\\ \\
We  
order the positive roots of $\fg$ as in Theorem \ref{bsb3}, 
thus for each $\pi \in {\overline{\bf P}}(\gamma)$   the odd root vector is the last factor of $e_{-\pi}$, that is we have $e_{-\pi}\in U(\fn^-_0)\fn^-_1$.
\bl \label{1car}If If $\gc$ is $\gl$-minimal, then $e_{-\gc} v_\gl$ occurs with non-zero coefficient in $\gth_\gc v_\gl$.
\el
\bpf Assume $\gc = \gep_r
- \gd_{s}.$ Then if $\ga = \gep_r - \gep_i$ with $r<i$,  or $\ga = \gd_j - \gd_s$ with $j<s$ we have $(\gl + \gr, \ga^\vee) \neq 0$, since $\gc$ is $\gl$-minimal.  Thus by Equation \eqref{98t}  and Theorem \ref{bsb3}, the  coefficient of $e_{-\gc} v_\gl$ 
 in $\gth_\gc v_\gl$ is  non-zero.  \epf
\bt \label{1boy} The isotropic root $\gc$ is independent at  $\gl$ if and only if $\gc$ is $\gl$-{minimal}. \et
\bpf Set $B=B(\gl)^{-\gc} $. If
$\gc$ is not $\gl$-{minimal}
then $\gc$ is not independent at  $\gl$
by Proposition \ref{1pot}. Suppose that $\gc$ is $\gl$-{minimal}
and
$$\theta_{\gc}v_{\gl}\in \sum_{\gc' \in B}U(\fg)\theta_{\gc'} v_{\gl} = \sum_{\gc' \in B}U(\fn^-)\theta_{\gc'} v_{\gl}.$$ 
Hence by 
Lemma \ref{1car}, there are $a_\pi \in \ttk$ such  that 
$$e_{-\gc}v_{\gl}= \sum_{\pi\in {\overline{\bf P}}(\gc): e_{-\pi}  \neq e_{-\gc}}a_\pi e_{-\pi}v_{\gl}.$$
This contradicts the fact that the set $ \{e_{-\pi}v_{\gl}|{\pi\in {\overline{\bf P}}(\gc)}\}$ is a basis for $M(\gl)^{\gl-\gc}$.
\epf

\subsection{Survival of \v{S}apovalov elements in Kac modules.}\label{srv}
We have $\mathfrak{g}_1 = \mathfrak{g}_1^+ \oplus
\mathfrak{g}_1^-,$ where $\mathfrak{g}_1^+ $ (resp. $\mathfrak{g}_1^-$) is the set of block  upper (resp. lower) triangular matrices.
Let $\fh$ be the Cartan subalgebra of $\fg$ consisting of diagonal matrices, and set $\mathfrak{p} =  \mathfrak{g}_0  \oplus \mathfrak{g}_1^{+}.$ Next let
 \begin{eqnarray*}
  P^+ &=&
   \{ \lambda \in \mathfrak{h}^*|(\lambda, \alpha^\vee) \in \mathbb{Z}, (\lambda, \alpha^\vee) \geq 0 \quad
 \mbox{for all}\quad \alpha \in \Delta^+_0 \}
 \end{eqnarray*}
For $\gl \in P^+,$ let $L^0(\lambda)$ be the (finite dimensional) simple  $\mathfrak{g}_0$-module with highest weight $\gl$.  Then
$L^0(\lambda)$ is naturally a $\fp$-module and we define the {\it Kac module} $K(\gl)$ by
\[K(\lambda) =  U(\mathfrak{g}) \otimes_{U(\fp)}
L^0(\lambda).\]
By a corollary to the PBW Theorem, see \cite{M} Corollary 6.1.5,
we have as a vector space
\[K(\lambda) =  \Lambda(\mathfrak{g}^-_1) \otimes
L^0(\lambda).\]
The next result is well-known.  Indeed two
methods of  proof are given in  Theorem 4.37 of \cite{Br}.  The
second of these is based on Theorem 5.5 in \cite{S2}.  We give a short
proof using  Theorem  \ref{bab3}. Assume that the roots are ordered as in  Theorem \ref{bab3}, that is  with the odd root vector first.
\bt \label{1shapel} Let $\gamma = \epsilon_r -
\delta_s$, and suppose  $\lambda$ and $\lambda - \gc$
belong to $P^+$ and $(\lambda+\rho,\gc) = 0$. Let 
the highest weight vector $v_\gl$ be the highest weight vector in $K(\gl)$ of weight $\gl$.
Then \bi
\itema There is a non-zero map of $\fg$-modules $K(\gl-\gc) \lra K(\gl)$, sending 
$xv_{\gl-\gc}$ to $x\theta_\gamma v_\gl$.
\itemb
$$[K(\lambda):L(\lambda-\gc)] \neq 0.$$\ei
 \et   
 \bpf 
 Let $\theta_\gamma$ be the in the \v Sapovalov element for the root $\gc$.  
Because of the way the roots are ordered, applying the right side of \eqref{22s} to the highest weight vector $v_\gl\in K(\gl)$
 yields 
\[ \theta_\gamma v_\lambda \in  \fg^-_1 \otimes L^0(\lambda) \subseteq \Lambda \fg^-_1 \otimes L^0(\lambda). \]
Therefore to prove (a) it suffices to
show that the coefficient of $e_{-\gamma} v_\lambda= e_{-\gamma}\ot  v_\lambda$ in $\theta_\gamma v_\lambda$  is
nonzero.  
Now the term $e_{-\gamma}$ arises in the sum \eqref{22s} when $I\in \I$ is as small as possible, that is $I= \{m+s, r\}.$  
In this case we have, by \eqref{67t}  and \eqref{98st} 
\[ H_I=\prod_{k=1}^{m-r} (h_{\gep_{r} - \gep_{r+k}}+\gr(h_{\gep_{r} - \gep_{r+k}}) -1)
\prod^{s-1}_{j=1} (h_{\gd_{j} -\gd_{s}}   +  \rho(h_{\gd_{j} -\gd_{s}})+1).    \]
Thus  
the coefficient of $e_{-\gamma}v_\lambda$ is
\[  H_I(\gl) =\prod_{k=1}^{m-r} ((\lambda + \rho, \gep_{r} - \gep_{r+k})-1)
\prod^{s-1}_{j=1} ( (\lambda + \rho, \gd_{j} -\gd_{s}) +1).    \]
Since $\lambda \in P^+$, the only factors in this product that could be zero   arise when $k=1$ or $j=s-1$. 
Let $\gs =\gep_{r} - \gep_{r+1}$ and $\gt= \gd_{s-1} -\gd_{s}$. Then 
the only factors that could be zero are 
\[    (\lambda + \rho, \gs)-1 = (\lambda , \gs)
\mbox{ and  }  (\lambda + \rho, \gt) +1 = (\lambda, \gt).    \]
However $(\gc, \gs^\vee) = (\gc, \gt^\vee) =1$. So if this happens, then $(\gl-\gc, \gs^\vee) =-1$ or $ (\gl-\gc, \gt^\vee) =-1$, and $\lambda - \gamma \notin
P^+$. The claim in (b) follows immediately from (a).\epf

\end{document}